\documentclass[reqno,twoside,11pt]{amsart}
 
\hoffset - 1.5cm

\usepackage{graphicx}
\usepackage{amstext}
\usepackage[T1]{fontenc}
\usepackage{pstricks,pst-node,pst-text,pst-3d}
\usepackage{color}
\usepackage{amsmath}
\usepackage{amsthm}
\usepackage{amssymb}
\usepackage{eufrak}
\usepackage[mathscr]{euscript}
\usepackage{bbm}

\newtheorem{Theorem}{Theorem}[section]
\newtheorem{Fact}{Fact}[section]
\newtheorem{Lemma}{Lemma}[section]
\newtheorem{Proposition}{Proposition}[section]
\newtheorem{Corollary}{Corollary}[section]
\theoremstyle{definition}
\newtheorem{Definition}{Definition}[section]
\newtheorem{Example}{Example}[section]
\newtheorem{Remark}{Remark}[section]

\newcommand{\numsec}{\setcounter{Theorem}{0}\setcounter{Definition}{0}
\setcounter{Remark}{0} \setcounter{Lemma}{0} \setcounter{Fact}{0}
\setcounter{Proposition}{0} \setcounter{Corollary}{0}
\setcounter{Example}{0} \setcounter{equation}{0}
\setcounter{Property}{0}\renewcommand\theequation{\arabic{section}.\arabic{equation}}
\renewcommand\theTheorem{\arabic{section}.\arabic{Theorem}}
\renewcommand\theDefinition{\arabic{section}.\arabic{Definition}}
\renewcommand\theRemark{\arabic{section}.\arabic{Remark}}
\renewcommand\theLemma{\arabic{section}.\arabic{Lemma}}
\renewcommand\theFact{\arabic{section}.\arabic{Fact}}
\renewcommand\theProposition{\arabic{section}.\arabic{Proposition}}
\renewcommand\theCorollary{\arabic{section}.\arabic{Corollary}}
\renewcommand\theExample{\arabic{section}.\arabic{Example}}
\renewcommand\theProperty{\arabic{section}.\arabic{Property}}}
\newcommand{\numsubsec}{\setcounter{Theorem}{0}\setcounter{Definition}{0}
\setcounter{Remark}{0} \setcounter{Lemma}{0} \setcounter{Fact}{0}
\setcounter{Proposition}{0} \setcounter{Corollary}{0}
\setcounter{Example}{0} \setcounter{equation}{0}
\setcounter{Property}{0}\renewcommand\theequation{\arabic{section}.\arabic{subsection}.\arabic{equation}}
\renewcommand\theTheorem{\arabic{section}.\arabic{subsection}.\arabic{Theorem}}
\renewcommand\theDefinition{\arabic{section}.\arabic{subsection}.\arabic{Definition}}
\renewcommand\theRemark{\arabic{section}.\arabic{subsection}.\arabic{Remark}}
\renewcommand\theLemma{\arabic{section}.\arabic{subsection}.\arabic{Lemma}}
\renewcommand\theFact{\arabic{section}.\arabic{subsection}.\arabic{Fact}}
\renewcommand\theProposition{\arabic{section}.\arabic{subsection}.\arabic{Proposition}}
\renewcommand\theCorollary{\arabic{section}.\arabic{subsection}.\arabic{Corollary}}
\renewcommand\theExample{\arabic{section}.\arabic{subsection}.\arabic{Example}}
\renewcommand\theProperty{\arabic{section}.\arabic{subsection}.\arabic{Property}}}

\newcommand{\numsubsubsec}{\setcounter{Theorem}{0}\setcounter{Definition}{0}
\setcounter{Remark}{0} \setcounter{Lemma}{0} \setcounter{Fact}{0}
\setcounter{Proposition}{0} \setcounter{Corollary}{0}
\setcounter{Example}{0} \setcounter{equation}{0}
\setcounter{Property}{0}\renewcommand\theequation{\arabic{section}.\arabic{subsection}.\arabic{subsubsection}.\arabic{equation}}
\renewcommand\theTheorem{\arabic{section}.\arabic{subsection}.\arabic{subsubsection}.\arabic{Theorem}}
\renewcommand\theDefinition{\arabic{section}.\arabic{subsection}.\arabic{subsubsection}.\arabic{Definition}}
\renewcommand\theRemark{\arabic{section}.\arabic{subsection}.\arabic{subsubsection}.\arabic{Remark}}
\renewcommand\theLemma{\arabic{section}.\arabic{subsection}.\arabic{subsubsection}.\arabic{Lemma}}
\renewcommand\theFact{\arabic{section}.\arabic{subsection}.\arabic{subsubsection}.\arabic{Fact}}
\renewcommand\theProposition{\arabic{section}.\arabic{subsection}.\arabic{subsubsection}.\arabic{Proposition}}
\renewcommand\theCorollary{\arabic{section}.\arabic{subsection}.\arabic{subsubsection}.\arabic{Corollary}}
\renewcommand\theExample{\arabic{section}.\arabic{subsection}.\arabic{subsubsection}.\arabic{Example}}
\renewcommand\theProperty{\arabic{section}.\arabic{subsection}.\arabic{subsubsection}.\arabic{Property}}}

\newcommand{\ba}{\begin{array}}
\newcommand{\bc}{\begin{center}}
\newcommand{\bd}{\begin{description}}
\newcommand{\bdm}{\begin{displaymath}}
\newcommand{\be}{\begin{enumerate}}
\newcommand{\beq}{\begin{equation}}
\newcommand{\bdf}{\begin{Definition}}
\newcommand{\bex}{\begin{Example}}
\newcommand{\bft}{\begin{Fact}}
\newcommand{\bl}{\begin{Lemma}}
\newcommand{\bp}{\begin{Proposition}}
\newcommand{\br}{\begin{Remark}}
\newcommand{\bt}{\begin{Theorem}}
\newcommand{\bco}{\begin{Corollary}}
\newcommand{\bh}{\begin{Hipothesis}}
\newcommand{\ea}{\end{array}}
\newcommand{\ec}{\end{center}}
\newcommand{\ed}{\end{description}}
\newcommand{\edm}{\end{displaymath}}
\newcommand{\ee}{\end{enumerate}}
\newcommand{\eeq}{\end{equation}}
\newcommand{\edf}{\end{Definition}}
\newcommand{\eex}{\end{Example}}
\newcommand{\eft}{\end{Fact}}
\newcommand{\el}{\end{Lemma}}
\newcommand{\ep}{\end{Proposition}}
\newcommand{\er}{\end{Remark}}
\newcommand{\et}{\end{Theorem}}
\newcommand{\eco}{\end{Corollary}}
\newcommand{\eh}{\end{Hipothesis}}

\newcommand{\bH}{\mathbb{H}}
\newcommand{\bI}{\mathbb{I}}

\newcommand{\bN}{\mathbb{N}}

\newcommand{\bR}{\mathbb{R}}

\newcommand{\bU}{\mathbb{U}}
\newcommand{\bV}{\mathbb{V}}
\newcommand{\bW}{\mathbb{W}}
\newcommand{\bX}{\mathbb{X}}

\newcommand{\bZ}{\mathbb{Z}}

\newcommand{\cC}{\mathcal{C}}

\newcommand{\cF}{\mathcal{F}}

\newcommand{\cH}{\mathcal{H}}
\newcommand{\cI}{\mathcal{I}}

\newcommand{\cN}{\mathcal{N}}

\newcommand{\cT}{\mathcal{T}}

\newcommand{\im}{\mathrm{ im \;}}
\newcommand{\diag}{\mathrm{ diag \;}}
\numberwithin{equation}{section} \errorcontextlines=0

\newcommand{\vp}{\varphi}

\newcommand{\cl}{\mathrm{cl}}

\newcommand{\sone}{S^1}
\newcommand{\sotwo}{SO(2)}

\newcommand{\sot}{SO(3)}
\newcommand{\ds}{\displaystyle}
\newcommand{\nt}{\noindent}

\newcommand{\subg}{\overline{\mathrm{sub}}(G)}
\newcommand{\subgc}{\overline{\mathrm{sub}}[G]}
\newcommand{\subh}{\overline{\mathrm{sub}}(H)}
\newcommand{\subhc}{\overline{\mathrm{sub}}[H]}

\newcommand{\sub}{\overline{\mathrm{sub}}}

\newcommand{\CW}{\mathrm{CW}}
\newcommand{\gcw}{G\text{-CW}}

\newcommand{\chig}{\chi_G}
\newcommand{\chih}{\chi_H}

\newcommand{\upsh}{\Upsilon_H}

\newcommand{\ci}{\cC\cI}
\newcommand{\cig}{\ci_G}
\newcommand{\cih}{\ci_H}

\newcommand{\sci}{\mathscr{CI}}
\newcommand{\dgh}{\nabla_H\text{-}\mathrm{deg}}

\newcommand{\lpm}{\lambda_{\pm}}

\newcommand{\h}{\mathbb{H}}

\usepackage{todonotes}
\usepackage{setspace}

\begin{document}

\title[Periodic orbits]{Symmetric Liapunov center theorem \\ for orbit with nontrivial isotropy group}

\author{Marta Kowalczyk$^{1)}$}
\address{$^{1), 2)}$Faculty of Mathematics and Computer Science \\ Nicolaus Copernicus University in Toru\'n\\
PL-87-100 Toru\'{n} \\ ul. Chopina $12 \slash 18$ \\
Poland}

\author{Ernesto P\'erez-Chavela}
\address{Departamento de Matem\'aticas, Instituto Tecnol\'ogico Aut\'onomo de M\'exico (ITAM),
R\'io Hondo 1, Col. Progreso Tizap\'an, 01080 Mexico, DF, Mexico}

\author{S{\l}awomir Rybicki$^{2)}$}
%\address{Faculty of Mathematics and Computer Science\\
%Nicolaus Copernicus University \\
%PL-87-100 Toru\'{n} \\ ul. Chopina $12 \slash 18$ \\
%Poland}

\email{martusia@mat.umk.pl (M. Kowalczyk)}
\email{ernesto.perez@itam.mx (E. P\'erez-Chavela)}
\email{rybicki@mat.umk.pl (S. Rybicki)}

\date{\today}


\begin{abstract}
In this article we prove two versions of the Liapunov center theorem for symmetric potentials. We consider a~second order autonomous system $\ddot q(t)=-\nabla U(q(t))$ in the presence of symmetries of a compact Lie group $\Gamma$ acting linearly on $\bR^n.$ We look for non-stationary periodic solutions of this system in a~neighborhood of an orbit of critical points of the potential $U.$  Our results generalize that of \cite{PCHRYST, PCHRYST2}. As a topological tool we use an infinite-dimensional generalization of the equivariant Conley index due to Izydorek, see \cite{IZYDOREK}.
\end{abstract}
%%%%%%%%%%%%%%%%%%%%%%%%%%%%%%%%%%%%%%%%%%%%%%%%%%%%%%%%%%%%%%%%%%%%%%%%%%%%%%%%%%%%%%%%%%%%%%%%%%%%%%%%%%%%
%%%%%%%%%%%%%%%%%%%%%%%%%%%%%%%%%%%%%%%%%%%%%%%%%%%%%%%%%%%%%%%%%%%%%%%%%%%%%%%%%%%%%%%%%%%%%%%%%%%%%%%%%%%%
%%%%%%%%%%%%%%%%%%%%%%%%%%%%%%%%%%%%%%%%%%%%%%%%%%%%%%%%%%%%%%%%%%%%%%%%%%%%%%%%%%%%%%%%%%%%%%%%%%%%%%%%%%%%
%%%%%%%%%%%%%%%%%%%%%%%%%%%%%%%%%%%%%%%%%%%%%%%%%%%%%%%%%%%%%%%%%%%%%%%%%%%%%%%%%%%%%%%%%%%%%%%%%%%%%%%%%%%%
%%%%%%%%%%%%%%%%%%%%%%%%%%%%%%%%%%%%%%%%%%%%%%%%%%%%%%%%%%%%%%%%%%%%%%%%%%%%%%%%%%%%%%%%%%%%%%%%%%%%%%%%%%%%
\maketitle
\section{introduction}\label{introduction}
\numsec
  
\nt
The Liapunov center theorem is one of the most significant theorems regarding the existence of periodic solutions of  ordinary differential equations in a neighborhood of stationary ones.
Consider a~second order autonomous system of the following form
\beq\label{eq0}
  \ddot q(t)=-\nabla U(q(t)),
\eeq
where the potential $U:\bR^n\rightarrow\bR$ is  of class $C^2$ and $0\in\bR^n$ is a non-degenerate critical point of $U$ which is not a local maximum, i.e. $\nabla U(0)=0,$ $\det\nabla^2U(0)\neq 0$ and  
$\sigma(\nabla^2 U(0)) \cap (0,+\infty) = \{\beta_1^2,\ldots, \beta_m^2\}$ for some $m \geq 1.$ 
Without loss of generality we assume that  $\beta_1>\beta_2>\ldots>\beta_m>0.$ Now suppose that there exists $\beta_{j_0}$ such that $\beta_j \slash \beta_{j_0}\not \in \bN$ for all $j=1,\ldots,j_0-1.$
Then the famous Liapunov center theorem states that there exists a~sequence $(q_k(t))$ of periodic solutions of the system \eqref{eq0} with amplitude tending to $0$ and a sequence $(T_k)$ of minimal periods such that $ T_k \rightarrow 2\pi/\beta_{j_0}$ as $ k\rightarrow+\infty,$  see for instance \cite{BERGER1, BERGER2}, \cite{HENRARD} and \cite{MAWI}.

A discussion of some generalizations of the  Liapunov center theorem one can find in \cite{PCHRYST}.

The goal of  this paper is to prove  the Liapunov center theorem in the presence of symmetries of the potential $U.$  Therefore, from now on, we  discuss  symmetric versions of the Liapunov center theorem.

Let  $\Omega \subset \bR^n$ be an open and $\Gamma$-invariant subset of $\bR^n$ where $\bR^n$ is considered as an orthogonal representation of a compact Lie group $\Gamma$. Assume that $q_0 \in \Omega$ is a critical point of the $\Gamma$-invariant potential $U : \Omega \to \bR$  of class $C^2$. Since for all $\gamma \in \Gamma$ the equality $U(\gamma  q_0)=U(q_0)$ holds and $\nabla U(q_0)=0$, the orbit 
$\Gamma(q_0)=\{\gamma q_0 : \gamma \in \Gamma\}$ consists of critical points of $U,$ i.e. $\Gamma(q_0) \subset (\nabla U)^{-1}(0).$ It is easy to see that $\dim \ker \nabla^2 U(q_0) \geq \dim \Gamma(q_0)  .$ The orbit $\Gamma(q_0)$  is called non-degenerate if 
$\dim \ker \nabla^2 U(q_0) = \dim \Gamma(q_0)  .$

For $\varepsilon>0$  by $\Gamma(q_0)_\varepsilon$ we understand an $\varepsilon$-neighborhood of the orbit $\Gamma(q_0),$ i.e. $\ds \Gamma(q_0)_\varepsilon =\bigcup_{q \in \Gamma(q_0)} B_\varepsilon(\bR^n,q),$ where $B_\varepsilon(\bR^n,q)$ denotes the open $\varepsilon$-ball centered at $q$ in $\bR^n.$

We are interested in finding non-stationary periodic solutions of the system \eqref{eq0} in a neighborhood of the orbit $\Gamma(q_0)$  of stationary solutions.
Note that if $\dim \Gamma \geq 1$ then it can happen that $\dim \Gamma(q_0) \geq 1,$ i.e. the critical point $q_0$ is not isolated in $(\nabla U)^{-1}(0)$. That is why for higher-dimensional orbits $\Gamma(q_0)$ one can not apply the classical Liapunov center theorem.

 In \cite{PCHRYST} we have proved the symmetric Liapunov center theorem for a non-degenerate orbit $\Gamma(q_0)$ of critical points of $U.$ More precisely, with the additional assumption that the isotropy group  $\Gamma_{q_0}=\{\gamma \in \Gamma : \gamma q_0=q_0\}$ is trivial   and that there is at least one positive eigenvalue of the Hessian  $\nabla^2 U(q_0),$ we have proved the existence of non-stationary periodic solutions of the system \eqref{eq0} in any neighborhood of the orbit $\Gamma(q_0).$ Moreover, we have controlled the minimal periods of these solutions in terms of positive eigenvalues of $\nabla^2 U(q_0),$ see Theorem 1.1 of  \cite{PCHRYST}.

In \cite{PCHRYST2} we have proved the symmetric Liapunov center theorem for a  minimal orbit $\Gamma(q_0).$ We have assumed that $\Gamma(q_0)$ is an isolated orbit of critical points of $U$ which is also an  orbit of minima of $U$ and that the isotropy group  $\Gamma_{q_0}$ is trivial. Requiring that there is at least one positive eigenvalue of  $\nabla^2 U(q_0)$ we have proved the existence of non-stationary periodic solutions of the system \eqref{eq0} in any neighborhood of the orbit $\Gamma(q_0).$ Moreover, we have controlled the minimal periods of these solutions in terms of positive eigenvalues of $\nabla^2 U(q_0),$ see Theorem 1.1 of  \cite{PCHRYST2}. We emphasize that  in this theorem  the orbit $\Gamma(q_0)$ can be degenerate, i.e. $\dim\ker\nabla^2U(q_0)>\dim\Gamma(q_0).$

Since the orbit $\Gamma(q_0)$ is $\Gamma$-homeomorphic to $\Gamma \slash \Gamma_{q_0}$ and in both theorems discussed above we have assumed that the isotropy group $\Gamma_{q_0}$ is trivial, the orbit $\Gamma(q_0)$ is $\Gamma$-homeomorphic to the group $\Gamma.$

As far as we know, there is no symmetric  Liapunov center theorem  for an orbit $\Gamma (q_0) $ of dimension at least $1$ with nontrivial isotropy group $\Gamma_{q_0}.$
Therefore the aim of this article is to prove two versions of such theorems.

Let $l\in\bN\cup\{0\}$ and by $T^l$ we understand the $l$-dimensional torus, i.e. \[T^l=\left\{\begin{array}{ll} \{e\}, & \text{if}\,\, l=0 \\[0.5em] \underbrace{S^1\times\cdots\times S^1}_{l\text{-times}}, & \text{if}\,\, l\neq 0 \end{array}\right..\] 
Let $H$ and $K$ be arbitrary groups. We write $H\approx K$ if the group $H$ is isomorphic to the group $K$.

Note that in both theorems formulated below if the group $\Gamma$ is abelian then the isotropy group $\Gamma_{q_0} $ can be arbitrary. On the other hand, if the group of symmetries $\Gamma$ is not abelian we assume that the isotropy group $\Gamma_{q_0}$ is isomorphic to a torus. A natural question arises whether these theorems can be strengthened by assuming that the isotropy group of $q_0$ is arbitrary. This question is at present far from being solved.

The following theorems significantly extend the class of potential applications. For example, 
if $\Gamma=\sot$ and the isotropy group $\Gamma_{q_0}$ is isomorphic to the circle group  $\sotwo \approx \sone,$ then the following spaces are homeomorphic:
$\Gamma(q_0),  \Gamma \slash \Gamma_{q_0}, \sot \slash \sotwo$ and $ S^2,$ i.e. the orbit $\Gamma(q_0)$ is homeomorphic to the two-dimensional sphere $S^2.$   We underline that this case is not covered by theorems proved in  \cite{PCHRYST,PCHRYST2}.

The theorem below is an extension of Theorem 1.1 of \cite{PCHRYST}, which we obtain assuming that  $l_0=0,$ i.e. $\Gamma_{q_0}=\{e\}.$

\bt \label{main-theo1} [Symmetric Liapunov center theorem for a non-degenerate orbit]
Let  $\Omega \subset \bR^n$ be an open and $\Gamma$-invariant subset of an orthogonal representation $\bR^n$ of a compact Lie group $\Gamma$.
Assume that $U:\Omega \to \bR$ is a $\Gamma$-invariant potential of class $C^2$ and $q_0 \in \Omega \cap (\nabla U)^{-1}(0)$. If moreover,
\be
\item $\Gamma$ is abelian or $\Gamma_{q_0}\approx T^{l_0}$ for some $l_0\in\bN\cup\{0\},$ 
\item $\dim\ker\nabla^2 U(q_0)=\dim \Gamma(q_0),$
\item $\sigma(\nabla^2 U(q_0)) \cap (0,+\infty) = \{\beta_1^2,\ldots, \beta_m^2\}$, $\beta_1>\beta_2>\ldots>\beta_m>0$ and $m\geq 1,$
\ee
then for any $\beta_{j_0}$ such that $\beta_j \slash \beta_{j_0} \not \in \bN $ for $j \neq j_0,$ there exists a sequence $(q_k(t))$ of periodic solutions of the system \eqref{eq0} with a sequence $(T_k)$ of minimal periods such that $T_k\to 2 \pi \slash \beta_{j_0}$ and   for any $\varepsilon>0$ there exists $k_0\in\bN$ such that $q_k([0,T_k])\subset \Gamma(q_0)_\varepsilon$ for all $k\geq k_0.$
\et

The following theorem  is a generalization of Theorem 1.1 of \cite{PCHRYST2}, which we obtain putting $l_0=0,$ i.e. $\Gamma_{q_0}=\{e\}.$

\bt\label{main-theo2} [Symmetric Liapunov center theorem for a minimal orbit]
Let  $\Omega \subset \bR^n$ be an open and $\Gamma$-invariant subset of an orthogonal representation $\bR^n$ of a compact Lie group $\Gamma$. 
Assume that $U:\Omega \to \bR$ is a $\Gamma$-invariant potential of class $C^2$ and $q_0 \in \Omega \cap (\nabla U)^{-1}(0)$. If moreover,
\be
\item $\Gamma$ is abelian or $\Gamma_{q_0}\approx T^{l_0}$ for some $l_0\in\bN\cup\{0\},$ 
\item $\Gamma(q_0)$ consists of minima of the potential $U$,
\item $\Gamma(q_0)$ is isolated in $(\nabla U)^{-1}(0)$,
\item $\sigma(\nabla^2 U(q_0)) \cap (0,+\infty) = \{\beta_1^2,\ldots, \beta_m^2\}$, $\beta_1>\beta_2>\ldots>\beta_m>0$ and $m\geq 1,$
\ee
then for any $\beta_{j_0}$ such that $\beta_j \slash \beta_{j_0} \not \in \bN $ for $j \neq j_0,$ there exists a sequence $(q_k(t))$ of periodic solutions of the system \eqref{eq0}  with a sequence $(T_k)$ of minimal periods such that $T_k\to 2 \pi \slash \beta_{j_0}$ and  for any $\varepsilon>0$ there exists $k_0 \in \bN$ such that $q_k([0,T_k])\subset \Gamma(q_0)_\varepsilon$ for all $k\geq k_0.$
\et

How do we prove these theorems? As in articles \cite{PCHRYST,PCHRYST2}, we consider periodic solutions of the system \eqref{eq0} as 
orbits of critical points of a $(\Gamma \times \sone)$-invariant functional defined on suitable chosen infinite-dimensional orthogonal representation $\bH^1_{2\pi}$ of  the group $\Gamma \times \sone.$ To prove our theorems we apply techniques of equivariant bifurcation theory. 
To be more precise, we prove a change of the $(\Gamma \times \sone)$-equivariant Conley index, see \cite{IZYDOREK}, along the family of trivial orbits $\Gamma(q_0) \times (0,+\infty) \subset \bH^1_{2\pi} \times (0, +\infty)$ which implies bifurcation of non-stationary periodic solutions of the system \eqref{eq0}.

Suppose that $\Gamma_{q_0} \approx T^{l_0}$ for some $l_0 \in \bN.$ Since the pair $(\Gamma \times \sone, \Gamma_{q_0} \times \sone)$ is not admissible, the homomorphism $i^{\star} : U(\Gamma_{q_0} \times \sone) \to U(\Gamma \times \sone)$ of the Euler rings induced by the inclusion homomorphism $i : \Gamma_{q_0} \times \sone \to \Gamma \times 
\sone$ is not an injection. Therefore,  proving the theorems formulated above, we must perform more subtle and advanced calculations in the Euler ring $U(\Gamma_{q_0} \times \sone)$
than those which were done  in articles \cite{PCHRYST,PCHRYST2}.

After introduction our paper is organized as follows. In Section 2 we introduce the equivariant setting and review some of standard facts on equivariant topology and representation theory of compact Lie groups. In Subsection \ref{groups} we recall the definition of the Euler ring $U(G)$ of a compact Lie group $G.$ Moreover, the properties of the equivariant Euler characteristic $\chig(\bX) \in U(G)$ of a~finite pointed $G$-$\CW$-complex  $\bX$ are also discussed in this subsection. In Subsection \ref{torus} we look more closely at the $H$-equivariant Euler characteristic $\chih(S^{\bV})\in U(H)$ of the $H$-$\CW$-complex $S^{\bV}$ where $\bV$ is an orthogonal representation of $H \approx T^{l_0}.$
Section \ref{proofmain-theo12} contains the proofs of our main results. Subsection \ref{VarSetting} is dedicated to introducing the variational setting of the problem. Namely, we study periodic solutions of the system \eqref{eq0} as orbits of critical points of a~$(\Gamma \times \sone)$-invariant functional defined on the Hilbert space $\bH^1_{2\pi}.$  Sections \ref{proofthm1} and \ref{proofthm2} are devoted to the proofs of Theorems \ref{main-theo1} and \ref{main-theo2}, respectively.

\newpage
\section{preliminary results}\label{premresults}
\numsec

\subsection{Groups and their representations }
\label{groups}
\numsubsec
\nt 
In this section for the convenience of the reader we repeat the relevant material from 
\cite{KBO} and \cite{DIECK} without proofs, thus making our exposition self-contained.

Let $G$ stand for a~compact Lie group and $\bV=(\bR^n,\varsigma)$ be a finite-dimensional, real, orthogonal representation of $G$, that is a~pair consisting of the space $\bR^n$ and a~continuous homomorphism $\varsigma: G\rightarrow O(n),$
where $O(n)$ denotes the group of orthogonal matrices.
We call $\bV$ trivial if $\varsigma(g)=Id_n$ for any $g\in G,$ where $Id_n$ is the identity matrix. 
The linear action of $G$ on $\bV$ is given by $G \times \bV \ni (g,v) \mapsto \varsigma(g) v \in \bV.$
To shorten notation, we continue to write $gv$ for $\varsigma(g)v$ and by $v\in\bV$ we understand $v\in\bR^n.$
A subset $\Omega\subset\mathbb{V}$ is called $G$-invariant if for any $g\in G$ and $v\in\Omega$ we have $gv\in\Omega.$ 
By an orthogonal subrepresentation of $\bV$ we understand a linear subspace $\bW\subset\bV$ which is also a $G$-invariant set. Additionally, define $\bV^G=\{v \in \bV : gv=v \:\: \forall g \in G\}.$

Two orthogonal representations of $G,$ $\bV=(\bR^n,\varsigma)$ and $\bV'=(\bR^n,\varsigma')$ are equivalent 
if there exists a $G$-equivariant, linear isomorphism $L : \bV \to \bV',$ i.e. the isomorphism $L$ satisfying $L(gv)=gL(v)$ for any $g \in G$ and $v \in \bV.$ 
We denote it briefly by $\bV \approx_G \bV'.$ 

Let $(\cdot,\cdot)$ and $\| \cdot \|$ denote the standard scalar product and the standard norm on $\bR^n,$ respectively. 
For an orthogonal subrepresentation $\bW\subset\bV$ we define the orthogonal complement $\bW^\bot$ of $\bW$ as $\bW^\bot=\{v\in\bV:(v,w)=0\, \forall w\in\bW\} \subset\bV.$

By the sum of two orthogonal representations of $G,$ $\bV_1=(\bR^{n_1},\varsigma_1)$ and $\bV_2=(\bR^{n_2},\varsigma_2)$ we understand the representation $\bV_1\oplus\bV_2,$ i.e. $(\bR^{n_1+n_2},\varsigma_1\oplus\varsigma_2)$ where the continuous homomorphism $\varsigma_1\oplus\varsigma_2: G\rightarrow O(n_1+n_2)$ is given by $(\varsigma_1\oplus\varsigma_2)(g)=\diag (\varsigma_1(g),\varsigma_2(g)),\, g\in G.$ 

Fix $v_0 \in \bV^G$ and define $B_\varepsilon(\bV,v_0) = \{v \in \bV : \lVert v - v_0 \rVert < \varepsilon\}, $ $D_\varepsilon(\bV,v_0) = \cl B_\varepsilon(\bV,v_0),$ $ S_\varepsilon(\bV,v_0) = \partial D_\varepsilon(\bV,v_0) $ and $S^{\bV}_{\varepsilon,v_0} = D_\varepsilon(\bV,v_0)/S_\varepsilon(\bV,v_0).$ 
Since $\bV$ is an orthogonal representation of $G$ and $v_0 \in \bV^G,$ the sets $B_\varepsilon(\bV,v_0),D_\varepsilon(\bV,v_0),\, S_\varepsilon(\bV,v_0)$ and $S^{\bV}_{\varepsilon,v_0}$ are $G$-invariant.
For simplicity of notation, we write $B_\varepsilon(\bV),D_\varepsilon(\bV),$ $S_\varepsilon(\bV)$ and $S^\bV_\varepsilon$ for $v_0=0$ and $B(\bV),D(\bV),$ $S(\bV)$ and $S^\bV$ for $v_0=0$ and $\varepsilon=1.$

Let $\subg$  denote the set of closed subgroups of $G.$ Two subgroups $H, K \in\subg$ are called conjugate in $G$ if there exists $g\in G$ such that $H=gKg^{-1}.$ Conjugacy is an equivalence relation and the conjugacy class of $H\in\subg$ is denoted by $(H)_G.$
Moreover, $\subgc$ denotes the set of the conjugacy classes of closed subgroups of $G.$

If $v\in\mathbb{V}$ then $G_v=\{g\in G:gv=v\}\in\subg$ is the isotropy group of $v$ and $G(v)=\{gv:g\in G\}\subset\bV$ is the $G$-orbit through $v.$ Notice that the isotropy groups of points on the same $G$-orbit are conjugate in  $G.$

Fix $k,l\in \mathbb{N}\cup\{\infty\}$ and an open, $G$-invariant subset $\Omega\subset\mathbb{V}.$ A map $\varphi:\Omega\rightarrow\mathbb{R}$ of class $C^{k}$ is said to be a $G$-invariant $C^k$-potential if 
$\varphi(gv)=\varphi(v)$ for any $g\in G$ and $v\in\Omega.$ 
The set of $G$-invariant $C^{k}$-potentials is denoted by $C^k_G(\Omega,\mathbb{R}).$
A map $\psi:\Omega\rightarrow\mathbb{V}$ of class $C^{l}$ is called a $G$-equivariant $C^l$-map if 
$\psi(gv)=g\psi(v)$
for any $g\in G$ and $v\in\Omega.$ 
The set of $G$-equivariant $C^{l}$-maps is denoted by $C^l_G(\Omega,\mathbb{V}).$
For any $G$-invariant $C^k$-potential $\varphi$ the gradient of $\varphi$ denoted by $\nabla\varphi$ is $G$-equivariant $C^{k-1}$-map. Similarly, we use the symbol $\nabla^2\varphi$ to denote the Hessian of $\varphi.$

Let us  recall the notion of an admissible pair which was introduced in \cite{PCHRYST}.
\bdf\label{admissible}
Let $H\in\subg.$ A~pair $(G,H)$ is said to be admissible if for all $K_1,\, K_2\in\subh$ the following implication holds true: if 
$(K_1)_H\neq(K_2)_H \text{ then } (K_1)_G\neq(K_2)_G.$
\edf

\br
Note that a~pair $(G,H)$ is admissible if for all $K_1,\, K_2\in\subh$ the following equivalence holds:
$(K_1)_H\neq(K_2)_H \text{ iff } (K_1)_G\neq(K_2)_G.$
\er

\br \label{remarkamissible}
Let $G$ be  abelian  and $H\in\subg.$ Then the pair $(G,H)$ is admissible. 
Indeed, note that for all $K_1, K_2 \in \subh$ we have $(K_1)_H=(K_2)_H \text{ iff } K_1=K_2 \text{ iff } (K_1)_G=(K_2)_G .$ 
Generally speaking, the class of admissible pairs is very restrictive.
\er

\bex \label{remarkamissible1}
If $\Gamma$ is an abelian compact Lie group, $\Gamma' \in \overline{\mathrm{sub}}(\Gamma),$ $ G=\Gamma\times S^1$ and $H=\Gamma'\times S^1$ then the pair $(G,H)$ is admissible.
\eex
Let $\cF_\ast(G)$ denote the set of finite, pointed $\gcw$-complexes, see  \cite{DIECK} for the definition of $\gcw$-complex. 
The $G$-homotopy type of $X\in\cF_\ast(G)$ is denoted by $[X]_G$ and $\cF_\ast[G]$ is the set of $G$-homotopy types of finite, pointed $\gcw$-complexes. 
Let $F$ be a~free abelian group generated by $\cF_\ast[G]$ and $N$ be a~subgroup of $F$ generated by elements $[A]_G-[X]_G+[X/A]_G$ where $A,X\in\cF_\ast(G)$ and $A\subset X.$
Define $U(G)=F/N$ and let $\chig(X)$ be the class of an element $[X]_G\in F$ in $U(G).$
If $X$ is a~finite $\gcw$-complex without base point we put $\chig(X)=\chig(X^+)$ where $X^+=X\sqcup\{\ast\}$ and $\ast$ is a~separate point added such that $g\ast=\ast$ for all $g\in G.$
For $(X,\ast_X),(Y,\ast_Y)\in\cF_\ast(G)$ put $X\vee Y=(X\times\{\ast_Y\}\cup\{\ast_X\}\times Y)/\{(\ast_X,\ast_Y)\}$ and $X\wedge Y=X\times Y/X\vee Y.$ Then $X\vee Y,\, X\wedge Y\in \cF_\ast(G).$
Since $[X]_G-[X\vee Y]_G+[Y]_G=[X]_G-[X\vee Y]_G+[X\vee Y/X]_G\in N,$ we have $\chig(X)+\chig(Y)=\chig(X\vee Y).$
Additionally, the assignment $(X,Y)\mapsto X\wedge Y$ induces a product $U(G)\times U(G)\rightarrow U(G)$ given by the formula
$\chig(X)\star\chig(Y)=\chig(X\wedge Y).$ 

The proof of the following theorem one can find in \cite{DIECK}.

\bt \label{chi_G(X)}
The group $(U(G),+)$ is the free abelian group with basis $\chig(G/H^+)$ for $(H)_G\in\subgc.$
Moreover, if $X\in\cF_\ast(G)$ and $\bigcup \limits_{k=0}^{p}\{(k,(H_{j,k})_G): j=1, \ldots, q(k)\}$ is a~type of the cell decomposition of  $X$ then
$
\chig(X) = \sum\limits_{(H)_G \in \subgc }\eta^G_{(H)_G}(X)\cdot \chig(G/H^+)\in U(G)
$
where $\eta^G_{(H)_G}(X)=\sum\limits_{k=0}^p(-1)^{k}\,\nu(X,k,(H)_G)\in\bZ$ and $\nu(X,k,(H)_G)$ is the number of cells of dimension $k$ and of orbit type $(H)_G$ in $X.$
\et
The triple $(U(G),+,\star)$ is a~commutative ring with unity $\bI_{U(G)}=\chig($ $G/G^+)$ and it is called the Euler ring of $G,$  see \cite{DIECK0,DIECK} for more properties of  $U(G).$

Let $H \in \subg$ and $Y$ be a~$H$-space, see \cite{DIECK} for the definition of $H$-space.
Now define an action of $H$ on the product $G\times Y$ by the formula
$(h,(g,y))\mapsto (gh^{-1},hy)$ and let $G\times_H Y$ denote the space of $H$-orbits of this action.
We denote the $H$-orbit through $(g,y)$ briefly by $[g,y].$
The space $G\times_H Y$ is a~$G$-space with the following action $(g',[g,y]) \mapsto [g'g, y].$
For a~pointed $H$-space $Y$ and $G^+=G\sqcup\{\ast\}$ where $\ast$ is a~separate point added such that $g\ast=\ast$ for all $g\in G$ we have $G^+\wedge Y=G^+\times Y \slash G^+ \vee Y =G\times Y \slash G\times\{\ast\}.$
The space $G^+\wedge Y$ is a~pointed $H$-space with the following action $(h,(g,y)) \mapsto (gh^{-1},hy)$ and by
$G^+\wedge_HY$ we denote the orbit space of this action. 
Similarly, we write the $H$-orbit through $(g,y)$ as $[g,y].$
Note that $G^+\wedge_HY$ is a~pointed $G$-space with an action induced by the assignment $(g',[g,y])\mapsto[g'g,y].$

The point of the following theorem is that it allows to express the $G$-equivariant Euler characteristic of a $\gcw$-complex $G^+\wedge_HY$ in terms of the $H$-equivariant Euler characteristic of the 
$H$-CW-complex $Y.$ 
The theorem below was proved in \cite{PCHRYST}, see Theorems 2.2 and 2.3 of \cite{PCHRYST}.
\bt \label{theo2.2/3ofPCHRYST}
Fix $H\in\subg$ and  $Y\in\cF_\ast(H).$ If $\chih(Y)=\sum\limits_{(K)_H\in\overline{sub}[H]}\eta^H_{(K)_H}(Y)\cdot\chih(H/K^+)$ 
then
\be
\item  $G^+\wedge_HY\in\cF_\ast(G),$

\item $\ds \chig(G^+ \wedge_H Y)  =  \sum_{(K)_H \in \subhc}\eta^H_{(K)_H}(Y) \cdot \chig(G/K^+) = \sum_{(K)_G \in \subgc }\eta^G_{(K)_G}(G^+\wedge_HY)\cdot \chig(G/K^+) \in U(G)$ and
$ \eta^G_{(K)_G}(G^+\wedge_H Y) = \sum \limits_{(K_1)_H \in \subhc,(K_1)_G = (K)_G} \eta^H_{(K_1)_H}(Y) \in \bZ,$

\item if moreover, the pair $(G,H)$ is admissible then 
$\eta^G_{(K)_G}(G^+\wedge_H Y) = \eta^H_{(K)_H}(Y)$ 
and that is why the map  $U(H) \ni \chih(Y) \to \chig\left(G^+ \wedge_HY\right) \in  U(G)$ is injective.
\ee
\et

\bco \label{remarkadmis}
If $G$ is  abelian  and $H \in \subg$ then by Remark \ref{remarkamissible} the pair $(G,H)$ is admissible and   the map  $U(H) \ni \chih(Y) \to \chig\left(G^+ \wedge_HY\right) \in  U(G)$ is injective.
\eco 
Fix an open, $G$-invariant subset $\Omega\subset\bV.$  Let $\varphi \in C^2_G(\Omega,\mathbb{R})$ and $G(q_0'),\, G(q_0'') \subset \Omega$ be non-degenerate critical $G$-orbits of $\varphi,$ i.e. $G(q_0^{\nu}) \subset (\nabla \varphi)^{-1}(0)\cap\Omega$ and $\dim \ker \nabla^2\varphi(q_0^{\nu}) = \dim G(q_0^{\nu})$ where $\nu \in \{',''\}.$ Additionally, assume that $G_{q_0'}=G_{q_0''}=H \in \subg.$ 

Fix $\nu \in \{',''\}$ and note that $T_{q_0^{\nu}} \bV = T_{q_0^{\nu}} G(q_0^{\nu}) \oplus T_{q_0^{\nu}}G(q_0^{\nu})^\bot$ where $T_{q_0^{\nu}}G(q_0^{\nu})^\bot$ is an orthogonal representation of $H,$ and  define $\phi^{\nu}=\varphi_{|T_{q_0^{\nu}}G(q_0^{\nu})^\bot} \in C^2_H(T_{q_0^{\nu}}G(q_0^{\nu})^\bot,\bR).$ 

In the following result we express the $G$-equivariant Conley index $\cig(G(q_0^{\nu}),-\nabla\varphi)$ of the non-degenerate critical $G$-orbit $G(q_0^{\nu})$  
 in terms of the $H$-equivariant Conley index $\cih(\{q_0^{\nu}\},-\nabla \phi^{\nu})$ of the non-degenerate critical point $q_0^{\nu}$ of the potential $\varphi$ restricted to the space orthogonal to this orbit.
So this theorem gains in interest if we realize that this relation allows us to distinguish the $G$-equivariant Conley indexes of non-degenerate orbits considering only the potential restricted to the spaces orthogonal to these orbits. The proof of the following theorem one can find in  \cite{PCHRYST}.

\bt \label{thmconley}
Under the above assumptions, if $\nu \in \{',''\}$ then 
\be
\item \label{thmconleya} $\cih(\{q_0^{\nu}\},-\nabla\phi^{\nu})\in\cF_\ast[H],$
\item \label{thmconleyb} $\cig(G(q_0^{\nu}),-\nabla\varphi)=G^+\wedge_H \cih(\{q_0^{\nu}\}, -\nabla\phi^{\nu}) \in \cF_\ast[G],$
\item \label{thmconleyc} if the pair $(G,H)$ is admissible and if $\chih(\cih(\{q_0'\},-\nabla\phi'))\neq\chih(\cih(\{q_0''\},-\nabla\phi''))\in U(H),$ then
$\cig(G(q_0'),-\nabla\varphi)\neq \cig(G(q_0''),-\nabla\varphi)\in\cF_\ast[G]$ and $\chig(\cig(G(q_0'),-\nabla\varphi))\neq\chig(\cig(G(q_0''),-\nabla\varphi))\in U(G).$
\ee
\et

\br\label{remarkringhom}
Let $H\in\subg.$ The standard homomorphism $i:H\rightarrow G$ induces the ring homomorphism $i^\star:U(G)\rightarrow U(H).$
\er

Rabinowitz proved that the Brouwer index of an isolated critical point which is  a local minimum of a potential of class $C^1$  equals $1\in\bZ,$ see Lemma 1.1 of \cite{RAB4}. 
The following lemma is an analogue of the Rabinowitz result for the class of equivariant gradient maps. 
Instead of the Brouwer degree we use the degree for equivariant gradient maps  $\dgh(\cdot,\cdot)\in U(H),$  see \cite{GEBA} for the definition and properties of this degree. In the proof of the lemma below we use the relation between the degree for equivariant gradient maps and the equivariant Conley index, see \cite{GEBA}, instead of the Poincar\'e-Hopf theorem used by Rabinowitz.

\bl \label{rabinowitz}
Let $\bV$ be an orthogonal representation of a compact Lie group $H$ and $\varphi\in C^1_H(\bV,\bR).$ Assume that $0\in\bV$ is an isolated critical point and local minimum of $\varphi.$ Then 
$\dgh(\nabla\varphi ,B_{\varepsilon}(\bV))=\bI_{U(H)}\in U(H).$
\el
\begin{proof}
There is no loss of generality in assuming that $\varphi(0)=0.$
Since $0 \in \bV$ is an isolated critical point of $H$-invariant potential $\varphi,$ $\{0\} \subset \bV$ is an isolated invariant set in the sense of the $H$-equivariant Conley index theory. Therefore the $H$-equivariant Conley index  $\cih(\{0\},-\nabla\varphi)$ is well defined. 
Since $0 \in \bV$ is an isolated critical  point of $\vp$, we can fix  $\varepsilon > 0$  such that $\nabla \varphi(v) \neq 0$ for $D_{\varepsilon}(\bV)\backslash \{0\}.$ Choose $\ds 0 < c < \min_{v \in S_\varepsilon(\bV)} \vp(v)$ a regular value of $\vp$ and define  $A_c= \vp^{-1}((- \infty,c]) \cap D_{\epsilon}(\bV).$
Since $\vp(0)=0,$ $0 \in A_c.$ Moreover, by the choice of the regular value $c$   we have $A_c \subset B_{\varepsilon}(\bV)$ and 
$\partial A_c=\vp^{-1}(c)\cap D_\varepsilon(\bV)$  is a manifold of codimension $1.$ Since $0 \in \bV$ is an isolated minimum of $\vp$, $A_c$ is contractible to a point  by using the negative gradient flow corresponding to $\varphi.$ Therefore  the pair $(A_c,\emptyset)$ is a $H$-index pair for the  isolated invariant set $\{0\}.$
Summing up, we obtain  $\chih(\cih(\{0\},-\nabla\varphi))=(-1)^0\chih(H/H^+)=\bI_{U(H)} \in U(H).$
Applying the  equality $\chih(\cih(\{0\},-\nabla\varphi)) = \dgh(\nabla \varphi , B_{\varepsilon}(\bV)),$ see \cite{GEBA}, we  complete the proof.
\end{proof}

The following lemma is analogous to Lemma \ref{rabinowitz}. But instead of minimum of the potential $\vp$  we consider its maximum.

\bl \label{rabinowitzco}
Let $\bV$ be an orthogonal representation of $H$ and $\varphi\in C^1_H(\bV,\bR).$ Assume that $0\in\bV$ is an isolated critical point and local maximum of $\varphi.$ Then
$\dgh(\nabla\varphi ,B_{\varepsilon}(\bV))=\chih(S^\bV)\in U(H).$
\el
\begin{proof}
Without loss of generality we can assume that $\varphi(0)=0.$ 
Like in the proof of Lemma \ref{rabinowitz} it follows that the $H$-equivariant Conley index  $\cih(\{0\},-\nabla\varphi)$ is well defined. 
Since $0 \in \bV$ is an isolated critical point of $\vp,$ we can choose $\varepsilon > 0$  such that $\nabla \varphi(v) \neq 0$ for $D_{\varepsilon}(\bV)\backslash \{0\}.$
Choose $\ds \max_{v \in S_\varepsilon(\bV)} \vp(v)< -c < 0$ a regular value of $\varphi$
and define $A_c= \varphi^{-1}([-c,+\infty)) \cap D_{\varepsilon}(\bV) .$
By the choice of the regular value $-c$ we obtain that $A_c \subset B_{\varepsilon}(\bV)$  and $0 \in A_c$.
Since  $-c$ is the regular value of $\vp,$ we have $\partial A_c = \vp^{-1}(-c)\cap D_\varepsilon(\bV)$ is a manifold of codimension $1.$ 
Moreover, since $0 \in \bV$ is a maximum of $\vp$, the negative gradient flow corresponding to $\varphi$ is directed outwards on $\partial A_c$. 
Additionally, $A_c$ is  contractible  by using the gradient flow corresponding to $\varphi.$ 
Hence the pair $(A_c,\partial A_c)$ is a $H$-index pair for the isolated invariant set $\{0\}.$  
It follows that 
$\chih(\cih(\{0\},-\nabla\varphi))=\chih(A_c/\partial A_c)=\chih(D(\bV)/S(\bV))=\chih(S^\bV)\in U(H).$
Applying the  equality $\chih(\cih(\{0\},-\nabla\varphi)) = \dgh(\nabla \varphi , B_{\varepsilon}(\bV)),$ see \cite{GEBA}, we  complete the proof.
\end{proof}

%%%%%%%%%%%%%%%%%%%%%%%%%%%%%%%%%%%%%%%%%%%%%%%%%%%%%%%%%%%%%%%%%%%%%%%%%%%%%%%%%%%%%%%%%%%%%%%%%%%%%%%%%%

%%%%%%%%%%%%%%%%%%%%%%%%%%%%%%%%%%%%%%%%%%%%%%%%%%%%%%%%%%%%%%%%%%%%%%%%%%%%%%%%%%%%%%%%%%%%%%%%%%%%%%%%%%%%
%%%%%%%%%%%%%%%%%%%%%%%%%%%%%%%%%%%%%%%%%%%%%%%%%%%%%%%%%%%%%%%%%%%%%%%%%%%%%%%%%%%%%%%%%%%%%%%%%%%%%%%%%%%%
%%%%%%%%%%%%%%%%%%%%%%%%%%%%%%%%%%%%%%%%%%%%%%%%%%%%%%%%%%%%%%%%%%%%%%%%%%%%%%%%%%%%%%%%%%%%%%%%%%%%%%%%%%%%
%%%%%%%%%%%%%%%%%%%%%%%%%%%%%%%%%%%%%%%%%%%%%%%%%%%%%%%%%%%%%%%%%%%%%%%%%%%%%%%%%%%%%%%%%%%%%%%%%%%%%%%%%%%%

\subsection{Torus}\label{torus}
\numsubsec
 
\nt
Let $l\in\bN\cup\{0\}$ and recall that 
\[T^l=\left\{\begin{array}{ll} \{e\}, & \text{if}\,\, l=0 \\[0.5em] \underbrace{S^1\times\cdots\times S^1}_{l\text{-times}}, & \text{if}\,\, l\neq 0 \end{array}\right..\] 
By $e^{i\phi}\in T^l$ we mean $e^{i\phi}=(e^{i\phi_1},\ldots,e^{i\phi_l})$ for some $\phi=(\phi_1,\ldots,\phi_l),\phi_i\in [0,2\pi),i=1,\ldots,l.$
Fix $m\in\bZ^l$ and define $H_m=\{e^{i\phi}\in T^l:e^{i( m,\phi )}=1\}=
\{e^{i\phi}\in T^l: ( m,\phi ) \in 2\pi\bZ\}\in\sub(T^l).$ Since $T^l$ is an abelian group, $(H_m)_{T^l}=H_{m}.$ Notice that $H_m=H_{m'}$ if and only if $m=\pm m'$ for every $m,m'\in\bZ^l.$

Fix $m\in\bZ^l\backslash\{0\}$ and let homomorphisms $\varsigma_m:T^l\rightarrow O(2)$ and $\varsigma_0:T^l\rightarrow O(1)$ be given by
\beq
\varsigma_m(e^{i\phi})=\left[\ba{rr}\cos ( m,\phi ) & -\sin ( m,\phi ) \\ \sin ( m,\phi ) & \cos ( m,\phi  ) \ea\right],\, \varsigma_0(e^{i\phi})=1.
\eeq 
Write $\varsigma_m^k=\diag(\underbrace{\varsigma_m,\ldots,\varsigma_m}_{k\text{-times}})$ and $\varsigma_0^k=\diag(\underbrace{\varsigma_0,\ldots,\varsigma_0}_{k\text{-times}})$ for some $k\in\bN.$

The following theorem gives the classification of irreducible representations of the torus $T^l,$ see  \cite{GLG}.

\bt The real, irreducible representations of the torus $T^l$ are the following:
\be
\item $\bR[1,m]=(\bR^2,\varsigma_m),\, m\in\bZ^l\backslash \{0\},$
\item $\bR[1,0]=(\bR,\varsigma_0).$
\ee
Representations $\bR[1,m]$ and $\bR[1,m']$ are equivalent if and only if $m=\pm m'.$ Moreover,
\[T^l_v=\left\{\ba{ll}H_m, & \text{if}\,\, v\in\bR[1,m]\backslash\{0\}\\ T^l, & \text{if}\,\, v=0\ea\right. .\]
\et

Let $G$ be a~compact Lie group and $H \in \sub(G).$ 

From now on we assume that there is $l_0 >0$ such that  $H\approx T^{l_0}$ and that $\bV$ is an orthogonal representation of $H$. Then there exist
$r,k_0\in\bN\cup\{0\},$ $k_1,\ldots,k_r\in\bN$ and $  m_1,\ldots,m_r\in\bZ^{l_0}\backslash\{0\}$ such that
$\bV\approx_{H}(\bR^n,\varsigma)=\bR[k_0,0]\oplus\bR[k_1,m_1] \oplus \cdots \oplus \bR[k_r,m_r]$
where $m_i\neq\pm m_j$ for $i\neq j,$
i.e. $\varsigma=\diag(\varsigma_0^{k_0},\varsigma_{m_1}^{k_1},\ldots,\varsigma_{m_r}^{k_r}).$

The theorem below yields partial information about the $H$-equivariant Euler characteristic of the $H$-CW-complex $S^{\bV}.$ 
The proof of this theorem one can find in \cite{GLG}.
\bt \label{G-R}
Under the above assumptions, the following equality holds
\begin{multline}
\chi_{H}(S^\bV)=(-1)^{k_0}\left(\chi_{H}(H/H^+)-\sum\limits_{i=1}^{r}k_i\cdot\chi_{H}(H/H_{m_i}^+)\right)+\\ +
\sum\limits_{(K)_{H}\in\{(\mathcal{H})_{H}\in \sub[H]:\dim\mathcal{H}\leq l_0-2\}}\eta^{H}_{(K)_{H}}(S^\bV)\cdot\chi_{H}(H/K^+)\in U(H).
\end{multline}
\et

Now we define a number $\mathfrak{S}(S^\bV)$ as the sum of absolute values of the coefficients of $\chih(S^\bV)$ assigned to
generators $\chih(H/K^+)$ of $U(H)$ such that $\dim K=l_0-1.$
Note that the only generators which satisfy this condition and have non-zero coefficients in $\chih(S^\bV)$ are of the form $\chi_{H}(H/H_{m_i}^+)$ for $i=1,\ldots,r.$

\bdf \label{sumcoeff}
Define $\mathfrak{S}(S^\bV)$ by
$\mathfrak{S}(S^\bV)=\sum\limits_{i=1}^r k_i\in\bN.$ It is understood that if $\bV$ is a trivial representation of $H$ we put $\mathfrak{S}(S^\bV)=0.$
\edf
Let $\bW$ be an orthogonal representation of $H$ equivalent to $\bR[k_0',0]\oplus\bR[k_1',m_1']\oplus\cdots\oplus\bR[k_s',m_s']$ for some $s,k_0'\in\bN\cup\{0\},$ $k_1',\ldots,k_s'\in\bN$ and $m_1',\ldots,m_s'\in\bZ^{l_0}\backslash\{0\}$
where $m_i'\neq\pm m_j'$ for $i\neq j.$
\br \label{remarkunique}
The decompositions of the orthogonal representations of $H,$  $\bV$ and $\bW$ are unique up to order of elements.
Note that $\bV\oplus\bW$ is equivalent to $\bR[k_0+k_0',0]\oplus\bR[k_1,m_1]\oplus\cdots\oplus\bR[k_r,m_r]\oplus\bR[k_1',m_1']\oplus\cdots\oplus\bR[k_s',m_s'].$
It is clear that there exist $t,k_0''\in\bN\cup\{0\},$ $k_1'',\ldots,k_t''\in\bN$ and $m_1'',\ldots,m_t''\in\bZ^{l_0}\backslash\{0\}$ such that
$\bV\oplus\bW\approx_{H}\bR[k_0'',0]\oplus\bR[k_1'',m_1'']\oplus\cdots\oplus\bR[k_t'',m_t''],$ $k_0''=k_0+k_0'$ and $\sum\limits_{i=1}^{t}k_i''=\sum\limits_{i=1}^{r}k_i+\sum\limits_{i=1}^{s}k_i'$
where $m_i''\neq\pm m_j''$ for $i\neq j.$
\er
\bl \label{propofsum}
Under the above assumptions, the following conditions hold.
\be
\item \label{propofsuma} The number $\mathfrak{S}(S^\bW)=0$ if and only if $\bW$ is a~trivial representation of $H,$ i.e. $\bW \approx_H \bR[k_0',0].$
\item \label{propofsumb} $\mathfrak{S}(S^{\bV\oplus\bW})=\mathfrak{S}(S^\bV)+\mathfrak{S}(S^\bW).$
\item \label{propofsumc} If $\bW$ is a nontrivial representation of $H$ then $\mathfrak{S}(S^\bV)\neq \mathfrak{S}(S^{\bV\oplus\bW}).$
\ee
\el
\begin{proof}
Condition (1) is obvious and condition (3) follows from (1) and (2). We only need to show (2). First of all, by Theorem \ref{G-R} and Remark \ref{remarkunique}, we have
\begin{align*}
\chi_{H}(S^{\bV\oplus\bW})=(-1)^{k_0''}\left(\chi_{H}(H/H^+)-\sum\limits_{i=1}^{t}k_i''\cdot\chi_{H}(H/H_{m_i''}^+)\right) +\\ +
\sum\limits_{(K)_{H}\in\{(\mathcal{H})_{H}\in \sub[H]:\dim\mathcal{H}\leq l_0-2\}}\eta^{H}_{(K)_{H}}(S^{\bV\oplus\bW})\cdot\chi_{H}(H/K^+)
\end{align*}
and $\sum\limits_{i=1}^{t}k_i''=\sum\limits_{i=1}^{r}k_i+\sum\limits_{i=1}^{s}k_i'.$
Then $\mathfrak{S}(S^{\bV\oplus\bW})=\sum\limits_{i=1}^{t}k_i''=\mathfrak{S}(S^\bV)+\mathfrak{S}(S^\bW),$ which is our claim.
\end{proof}
The following lemma can be easily deduced from Corollary 3.2 of \cite{GLG}.
\bl
The following conditions are equivalent:
\be
\item $\chih(S^\bV)=\chih(S^{\bV \oplus \bW}),$
\item $\bW$ is a trivial and even-dimensional representation of $H$.
\ee
\el
We can now formulate the relation between the number $\mathfrak{S}(S^\bV)$ and the $H$-equivariant Euler characteristic $\chih(S^\bV).$

\bl \label{col1}
If $\mathfrak{S}(S^\bV)\neq \mathfrak{S}(S^\bW)$ then $\chih(S^\bV)\neq \chih(S^\bW).$
\el
\begin{proof}
Without loss of generality we assume that $\mathfrak{S}(S^\bV)\neq 0.$  
Since $\mathfrak{S}(S^\bV)\neq 0,$ $\bV$ is a nontrivial representation of $H.$ Suppose, contrary to our claim, that $\mathfrak{S}(S^\bV)\neq \mathfrak{S}(S^\bW)$ and $\chih(S^\bV)=$ $ \chih(S^\bW).$ 
Since the representation  $\bV$ is nontrivial and $\chih(S^\bV)= \chih(S^\bW),$  it follows that $\bW$ is a~nontrivial representation of $H.$ 
By Theorem \ref{G-R}, we have
\begin{align*}
\chi_{H}(S^\bV)=& (-1)^{k_0}\left(\chi_{H}(H/H^+)-\sum\limits_{i=1}^{r}k_i\cdot\chi_{H}(H/H_{m_i}^+)\right)+\\ &
+\sum\limits_{(K)_{H}\in\{(\mathcal{H})_{H}\in \sub[H]:\dim\mathcal{H}\leq l_0-2\}}\eta^{H}_{(K)_{H}}(S^\bV)\cdot\chi_{H}(H/K^+),\\
\chi_{H}(S^\bW)= & (-1)^{k_0'}\left(\chi_{H}(H/H^+)-\sum\limits_{i=1}^{s}k_i'\cdot\chi_{H}(H/H_{m_i'}^+)\right)+\\ &
+\sum\limits_{(K)_{H}\in\{(\mathcal{H})_{H}\in \sub[H]:\dim\mathcal{H}\leq l_0-2\}}\eta^{H}_{(K)_{H}}(S^\bW)\cdot\chi_{H}(H/K^+).
\end{align*}
Therefore the equality $\chih(S^\bV)= \chih(S^\bW)$ implies that $k_0$ and $k_0'$ are of the same parity and 
\beq \label{eqcod1}
\sum \limits_{i=1}^{r} k_i \cdot \chi_{H}(H/H_{m_i}^+) = \sum \limits_{i=1}^{s}k_i' \cdot \chi_{H}(H/H_{m_i'}^+).
\eeq 
Note that the sets $\{\chi_{H}(H/H_{m'_1}^+), \ldots, \chi_{H}(H/H_{m'_s}^+)\}$ and $\{\chi_{H}(H/H_{m_1}^+), \ldots, \chi_{H}(H/H_{m_r}^+)\}$ consist of  linearly independent elements of $U(H).$
We claim that   $\{\chi_{H}(H/H_{m'_1}^+), \ldots, \chi_{H}(H/H_{m'_s}^+)\} = \{\chi_{H}(H/H_{m_1}^+), \ldots, \chi_{H}(H/H_{m_r}^+)\},$ $r=s$ and $\{k_1',\ldots,k_s'\} = \{k_1,\ldots,k_r\}.$ 
Suppose, contrary to our claim, that  there exists $\chi_{H}(H/H_{m_{i_0}'}^+) \not \in\{\chi_{H}(H/H_{m_1}^+), \ldots, \chi_{H}(H/H_{m_r}^+)\}.$ 
Hence $\ds \chi_{H}(H/H_{m_{i_0}'}^+) \notin \bigcup_{i=1,i \neq i_0}^s \{\chi_{H}(H/H_{m'_i}^+)\} \cup \bigcup_{i=1}^r \{\chi_{H}(H/H_{m_i}^+)\}.$
By the equality \eqref{eqcod1}, we obtain 
$k'_{i_0} \chi_{H}(H/H_{m_{i_0}'}^+) = \sum \limits_{i=1}^{r} k_i \cdot \chi_{H}(H/H_{m_i}^+) - \sum \limits_{i=1, \\ i \neq i_0}^{s}k_i' \cdot \chi_{H}(H/H_{m_i'}^+),$ i.e. the basis element $\chi_{H}(H/H_{m_{i_0}'}^+)$
is presented as a linear combination of basis elements,   see Theorem  \ref{chi_G(X)},
a contradiction. We have just proved that 
$$\{\chi_{H}(H/H_{m'_1}^+), \ldots, \chi_{H}(H/H_{m'_s}^+)\} \subset \{\chi_{H}(H/H_{m_1}^+), \ldots, \chi_{H}(H/H_{m_r}^+)\}.$$ 
The reverse conclusion is proved in the same way.
Since  $\{\chi_{H}(H/H_{m'_1}^+), \ldots, \chi_{H}(H/H_{m'_s}^+)\} = \{\chi_{H}(H/H_{m_1}^+), \ldots, \chi_{H}(H/H_{m_r}^+)\}$ and by the equality \eqref{eqcod1}, we have $r=s,$ for all $j=1,\ldots,r,$ $k_j =k'_j$ and $ m_j=\pm m'_j.$
 In consequence, $\mathfrak{S}(S^\bV)=\sum\limits_{i=1}^r k_i=\mathfrak{S}(S^\bW),$  a contradiction.
\end{proof}
As a direct consequence of Lemmas \ref{propofsum}.\eqref{propofsumc} and \ref{col1} we obtain the following corollary.

\bco \label{cor1}
If $\bW$ is a nontrivial representation of $H$ then $\chih(S^\bV)\neq\chih(S^{\bV\oplus\bW}).$
\eco
The crucial role in the following theorem plays the connection of
the $G$-equivariant Euler characteristic of the $G$-CW-complex $G^+\wedge_{H}S^{\bV}$ with the $H$-equivariant Euler characteristic of the $H$-CW-complex $S^\bV,$ see Theorem \ref{theo2.2/3ofPCHRYST}.
Notice that the pair $(G,H)$ does not need to be admissible, i.e. it can happen that there exist two subgroups $K_1,K_2\in\subg$ such that $(K_1)_G=(K_2)_G$ and $(K_1)_H\neq(K_2)_H.$

\bt \label{theo_rozne}
If $\bW$ is a~nontrivial representation of $H$ then the following inequality holds true $\chig(G^+\wedge_{H}S^{\bV})\neq\chig(G^+\wedge_{H}S^{\bV\oplus\bW}).$
\et
\begin{proof}
On the contrary, suppose that $\bW$ is a~nontrivial representation of $H$ and $\chig(G^+\wedge_{H}S^{\bV})=\chig(G^+\wedge_{H}S^{\bV\oplus\bW}).$ By Lemma \ref{propofsum}.\eqref{propofsumc}, we get $\mathfrak{S}(S^\bV)\neq \mathfrak{S}(S^{\bV\oplus\bW}).$ Note that the operation of conjunction preserves the dimensions of subgroups, i.e. $\dim K=\dim gKg^{-1}$ for every $g\in G$ and $K\in\subg.$
Thus Theorem \ref{theo2.2/3ofPCHRYST} now leads to

\begin{align*}
\chi_{G}(G^+\wedge_{H}S^{\bV})= &
(-1)^{k_0}\left(\chig(G/H^+)-\sum\limits_{i=1}^{r}k_i\cdot\chig(G/H_{m_i}^+)\right)+\\ &
+ \sum \limits_{(K)_G \in \{(\cH)_G \in \sub[G]: \dim \cH   \leq l_0 - 2 \}} \eta^{G}_{(K)_G}(G^+\wedge_{H}S^{\bV}) \cdot \chig(G/K^+),\\
\end{align*}

\begin{align*}
\chig(G^+\wedge_{H} S^{\bV \oplus \bW})= &
(-1)^{k_0+k_0'} \left (\chig(G/H^+) - \sum \limits_{i=1}^{r} k_i \cdot \chig(G/H_{m_i}^+)-\sum \limits_{j=1}^{s}k_j' \cdot \chig(G/H_{m_j'}^+)\right)+\\ &
+ \sum \limits_{(K)_G\in\{(\cH)_G\in \subgc: \dim \cH \leq l_0-2\}} \eta^{G}_{(K)_G}(G^+\wedge_{H}S^{\bV \oplus \bW}) \cdot \chig (G/K^+),
\end{align*}
and so $k_0'$ is an even number.
The equality $\chig(G^+\wedge_{H}S^{\bV})=\chig(G^+\wedge_{H}S^{\bV\oplus\bW})$ implies that
\[\sum\limits_{i=1}^{r}k_i\cdot\chig(G/H_{m_i}^+)=
\sum\limits_{i=1}^{r}k_i\cdot\chig(G/H_{m_i}^+)+\sum\limits_{j=1}^{s}k_j'\cdot\chig(G/H_{m_j'}^+),\] and consequently
$\sum\limits_{j=1}^{s}k_j'\cdot\chig(G/H_{m_j'}^+)=0\in U(G).$
By Theorem \ref{chi_G(X)}, we obtain  $k_j'=0$ for any $j,$ which gives $\mathfrak{S}(S^\bV)=\sum\limits_{i=1}^rk_i=\mathfrak{S}(S^{\bV\oplus\bW}),$ a contradiction.
\end{proof}

Recall that $\bV$ is an orthogonal representation of $T^{l_0}.$
Define $\bV_k = \{a\cos kt + b\sin kt: a,b\in\bV\}, k \geq 0.$ 
The action of the group $T^{l_0}\times S^1$ on $\bV_k$ is defined as follows
$$ (T^{l_0}\times \sone) \times \bV_k \ni ((e^{i\phi},e^{i\theta}),a\cos kt + b\sin kt) \mapsto  \varsigma(e^{i\phi})a\cos k(t+\theta) + \varsigma(e^{i\phi}) b \sin k (t+\theta).$$
Thus $\bV_k$ is an orthogonal representation of  $T^{l_0}\times \sone$.

\bl \label{lem_fixedpoints}
Under the above assumptions, $\bV_k^{T^{l_0}\times \sone}=\{0\}, k \geq 0.$
\el
\begin{proof}
Notice that $0 \in \bV_k^{T^{l_0} \times S^1}.$ Let $a \cos kt + b \sin kt \in \bV_k^{T^{l_0}\times S^1},$ i.e. for all $(e^{i\phi},e^{i\theta})\in T^{l_0}\times S^1$ we have
$(e^{i\phi},e^{i\theta})(a\cos kt + b\sin kt)=a\cos kt + b\sin kt.$
Therefore
\beq \label{eq1} \varsigma(e^{i\phi})a\cos k\theta + \varsigma(e^{i\phi})b\sin k\theta = a\,\, \text{and}\, -\varsigma(e^{i\phi})a\sin k\theta + \varsigma(e^{i\phi})b\cos k\theta = b.\eeq
Since the equalities \eqref{eq1} are satisfied for all $\phi$ and $\theta,$ for $\phi=0$ and $\theta=\pi/k$ we get $a=b=0,$ and the proof is complete.
\end{proof}

\section{Proofs of Theorems \ref{main-theo1} and \ref{main-theo2}}
\label{proofmain-theo12}
\numsec

To prove our main results we use techniques of equivariant bifurcation theory. Set $G=\Gamma \times \sone.$ We treat periodic solutions of the system \eqref{eq0} as  $G$-orbits of critical points of a $G$-invariant potential. As a topological tool we use the $G$-equivariant Conley index due to Izydorek, see  \cite{IZYDOREK}.  More precisely, we prove a change of the $G$-equivariant Conley index  along the trivial family $\Gamma(q_0) \times (0,+ \infty) \subset \bH^1_{2\pi} \times (0,+\infty).$ Such a change implies a local bifurcation of periodic solutions of the system $\ddot q(t)=-\nabla U(q(t)).$

\subsection{Variational setting}\label{VarSetting}
\numsubsec
In this section we introduce the variational setting for our problem, i.e. we consider periodic solutions of  the system \eqref{eq0} as critical $G$-orbits of $G$-invariant functionals. 
It is known that  there is a~one-to-one correspondence between $2\pi\lambda$-periodic solutions of the system \eqref{eq0} and $2\pi$-periodic solutions of the following family 
\beq \label{eqlambda}
\left\{\ba{rcl}
\ddot q(t) & = & -\lambda^2 \nabla U(q(t)) \\
      q(0) & = & q(2\pi) \\
 \dot q(0) & = & \dot q(2\pi)
\ea.
\right.
\eeq
Write
$
\h^1_{2\pi} = \{u : [0,2\pi] \rightarrow \bR^n : \text{ u is abs. continuous map, } u(0)=u(2\pi), \dot u \in L^2([0,2\pi],\bR^n)\}
$
and
$
\ds \langle u,v\rangle_{\h^1_{2\pi}} = \int_0^{2\pi} (\dot u(t), \dot v(t)) + (u(t),v(t)) \; dt.
$
Then $\left(\h^1_{2\pi},\langle\cdot,\cdot\rangle_{\h^1_{2\pi}}\right)$ is a separable Hilbert space which  is also an orthogonal representation of   $G$ with the following action
\[
G \times \h^1_{2\pi} \ni ((\gamma,e^{i \theta}),q(t)) \mapsto
\left\{\ba{rcl}
\gamma q(t+\theta), & \text{if}\,\, t+\theta< 2\pi \\ \gamma q(t+\theta - 2\pi), & \text{if}\,\, t+\theta\geq 2\pi
\ea\right. .
\] 
By $\|\cdot\|_{\bH^1_{2\pi}}$ we denote the norm induced by the inner product $\langle\cdot,\cdot\rangle_{\h^1_{2\pi}}.$
Let $\Phi : \h^1_{2\pi}\times (0,+\infty)  \to \bR$ be given by the formula
\beq  \label{functional}
\Phi(q,\lambda) = \int_0^{2\pi} \left( \frac{1}{2} \| \dot q(t) \|^2 -  \lambda^2U(q(t)) \right) \; dt,
\eeq
where $\lambda$ is treated as a~parameter.
Then $\Phi$ is a~$G$-invariant functional of class $C^2$ and solutions of the system \eqref{eqlambda}
correspond to solutions of the following equation
\beq\label{eqgrad}
\nabla_q\Phi(q,\lambda)=0.
\eeq
 Write $\h_0=\bR^n$ and $\h_k =\{a\cos kt+b\sin kt:a,b\in\bR^n\}$ for $k > 0.$ Note that
\beq \label{space}
\bH^1_{2\pi} = \overline{\h_0 \oplus \bigoplus_{k=1}^{\infty} \h_k}
\eeq
where the finite-dimensional space $\h_k$ is an orthogonal representation of $G$ for $k \geq 0.$

Since $q_0  \in \h^1_{2\pi}$ is a constant function,  $\mathcal{T}=G(q_0) \times(0,+\infty)= \Gamma(q_0)\times(0,+\infty) \subset \h^1_{2\pi}\times(0,+\infty).$ 
The family $\cT$ is called a~family of trivial solutions of the equation \eqref{eqgrad} while we call 
$\cN=\{(q,\lambda)\in  \h^1_{2\pi}\times(0,+\infty)\setminus\mathcal{T}:\nabla_q\Phi(q,\lambda)=0\} $ a set of nontrivial solutions.

Now we look for parameters which satisfy the necessary condition for the existence of local bifurcation, i.e. 
\beq \label{neccond}
\ker\nabla_q^2\Phi(q_0,\lambda)\cap \overline{\bigoplus_{k=1}^{\infty} \h_k}\neq \emptyset,
\eeq 
see Section 4 of \cite{PCHRYST} and Theorem 3.2.1 of \cite{PCHRYST2}.
The condition \eqref{neccond} is fulfilled if and only if 
$k^2-\lambda^2\beta_j^2=0$ for some $k \in \bN,$ see Lemma 5.1.1 of \cite{UGORY}. 
Let $\Lambda=\left\{\frac{k}{\beta_j}:k\in\bN,\, j=1,\ldots,m\right\}$ and then a~local bifurcation of solutions of the equation $\eqref{eqgrad}$ from the trivial family $\mathcal{T}$ can occur only from orbits $\Gamma(q_0)\times\Lambda\subset\mathcal{T}.$

Fix $\varepsilon > 0.$ By $G(q_0)_\varepsilon$ and $\Gamma(q_0)_\varepsilon$ we understand $\varepsilon$-neighborhoods of $G(q_0)$ and $\Gamma(q_0)$ in $\bH^1_{2\pi}$ and $\bH_0=\bR^n,$ respectively, i.e.  
$\ds G(q_0)_\varepsilon=\bigcup_{q \in G(q_0)}B_\varepsilon(\bH^1_{2\pi},q) \subset \bH^1_{2\pi}$ and $\ds \Gamma(q_0)_\varepsilon = \bigcup_{q \in \Gamma(q_0)} B_\varepsilon(\bH_0,q) \subset \bH_0=\bR^n.$

 It is known that a~change of the $G$-equivariant Conley index 
$\sci_{G}\left(G(q_0),-\nabla\Phi(\cdot,\lambda)\right)$
along the trivial family $\mathcal{T}$ implies the existence of a~local bifurcation of solutions of the equation \eqref{eqgrad}, where $\sci_{G}\left(G(q_0),-\nabla\Phi(\cdot,\lambda)\right)$
is an infinite-dimensional generalization of the $G$-equivariant Conley index which is a  $G$-homotopy type of a $G$-spectrum. This construction is due to Izydorek, see \cite{IZYDOREK} for more details.

Notice that the local bifurcations mentioned above are bifurcations in the function space $\bH^1_{2\pi}$ but more interesting phenomena are bifurcations in the phase space, and these kinds of bifurcations are the claims of our Theorems  \ref{main-theo1} and \ref{main-theo2}. The following lemma states that the existence of a local bifurcation in the phase space is a natural consequence of the occurrence of a local bifurcation in the function space.

\bl Under the above assumptions, if there exists a~local bifurcation of solutions of the equation \eqref{eqgrad} from the critical orbit $G(q_0),$ i.e. there exists a sequence $(q_k(t))$ of periodic solutions of the equation \eqref{eq0} such that for any $\varepsilon>0$ there exists $k_0\in\bN$ such that $q_k\in G(q_0)_\varepsilon$ for all $k\geq k_0,$ then
for any $\varepsilon>0$ there exists $k_0\in\bN$ such that $q_k([0,T_k])\subset \Gamma(q_0)_\varepsilon$ for all $k\geq k_0.$
\el 
\begin{proof}
Fix $\varepsilon > 0.$ Since there exists $k_0\in\bN$ such that $q_k\in G(q_0)_\varepsilon$ for all $k\geq k_0$ and $q_0$ is a constant function, we obtain $\gamma\in \Gamma$ such that 
$q_k\in B_\varepsilon(\bH^1_{2\pi},\gamma q _0),$ that is $\| q_k-\gamma q_0\|_{\h^1_{2\pi}} < \varepsilon.$ 
Therefore, by  Proposition 1.1 of \cite{MAWI}, for some $c>0$ we have 
$\ds \sup_{t\in[0,T_k]}\| q_k(t)- \gamma q_0\| =\|q_k- \gamma q_0\|_{\infty} \leq  c \| q_k- \gamma q_0\|_{\h^1_{2\pi}} < c\varepsilon.$   In consequence $q_k([0,T_k])\subset \Gamma(q_0)_{c\varepsilon}$ for all $k\geq k_0,$ which completes the proof.
\end{proof}
So to prove Theorems \ref{main-theo1} and \ref{main-theo2} we need to show a~change in the $G$-equivariant Conley index of the $G$-orbit $G(q_0).$

%%%%%%%%%%%%%%%%%%%%%%%%%%%%%%%%%%%%%%%%%%%%%%%%%%%%%%%%%%%%%%%%%%%%%%%%%%%%%%%%%%%%%%%%%%%%%%%%%%%%%%%%%%
%%%%%%%%%%%%%%%%%%%%%%%%%%%%%%%%%%%%%%%%%%%%%%%%%%%%%%%%%%%%%%%%%%%%%%%%%%%%%%%%%%%%%%%%%%%%%%%%%%%%%%%%%%%%
%%%%%%%%%%%%%%%%%%%%%%%%%%%%%%%%%%%%%%%%%%%%%%%%%%%%%%%%%%%%%%%%%%%%%%%%%%%%%%%%%%%%%%%%%%%%%%%%%%%%%%%%%%%%
%%%%%%%%%%%%%%%%%%%%%%%%%%%%%%%%%%%%%%%%%%%%%%%%%%%%%%%%%%%%%%%%%%%%%%%%%%%%%%%%%%%%%%%%%%%%%%%%%%%%%%%%%%%%
%%%%%%%%%%%%%%%%%%%%%%%%%%%%%%%%%%%%%%%%%%%%%%%%%%%%%%%%%%%%%%%%%%%%%%%%%%%%%%%%%%%%%%%%%%%%%%%%%%%%%%%%%%%%
%%%%%%%%%%%%%%%%%%%%%%%%%%%%%%%%%%%%%%%%%%%%%%%%%%%%%%%%%%%%%%%%%%%%%%%%%%%%%%%%%%%%%%%%%%%%%%%%%%%%%%%%%%%%
%%%%%%%%%%%%%%%%%%%%%%%%%%%%%%%%%%%%%%%%%%%%%%%%%%%%%%%%%%%%%%%%%%%%%%%%%%%%%%%%%%%%%%%%%%%%%%%%%%%%%%%%%%%%
%%%%%%%%%%%%%%%%%%%%%%%%%%%%%%%%%%%%%%%%%%%%%%%%%%%%%%%%%%%%%%%%%%%%%%%%%%%%%%%%%%%%%%%%%%%%%%%%%%%%%%%%%%%%
%%%%%%%%%%%%%%%%%%%%%%%%%%%%%%%%%%%%%%%%%%%%%%%%%%%%%%%%%%%%%%%%%%%%%%%%%%%%%%%%%%%%%%%%%%%%%%%%%%%%%%%%%%%%
%%%%%%%%%%%%%%%%%%%%%%%%%%%%%%%%%%%%%%%%%%%%%%%%%%%%%%%%%%%%%%%%%%%%%%%%%%%%%%%%%%%%%%%%%%%%%%%%%%%%%%%%%%%%

\subsection{Proof of Theorem \ref{main-theo1}} \label{proofthm1}
\numsubsec
We follow the notation used in \cite{PCHRYST}.
Fix $\beta_{j_0}$ satisfying the assumptions of Theorem \ref{main-theo1} and choose $\varepsilon>0$ such that $[\lambda_-,\lambda_+]\cap\Lambda=\{\frac{1}{\beta_{j_0}}\}$ where $\lpm=\frac{1\pm\varepsilon}{\beta_{j_0}}.$ To prove this theorem we have to show that $\sci_{G}\left(G(q_0),-\nabla\Phi(\cdot,\lambda_-)\right)\neq
\sci_{G}\left(G(q_0),-\nabla\Phi(\cdot,+)\right).$
The $G$-equivariant Conley index $\sci_{G}\left(G(q_0),-\nabla\Phi(\cdot,\lambda_\pm)\right)$ is a $G$-homotopy type of a $G$-spectrum $\left(\mathcal{E}_{n,\pm}\right)_{n=n_0}^{\infty}$
where $\mathcal{E}_{n,\pm}=\cig(G(q_0),-\nabla\Phi^n(\cdot,\lambda_\pm))$ and $\ds \Phi^n(\cdot,\lambda_\pm) = \Phi(\cdot,\lambda_\pm)_{| \oplus_{k=0}^{n} \bH_k} : \bigoplus_{k=0}^{n} \bH_k \rightarrow \bR.$
It follows that for any $n\geq n_0$ we have 
\begin{equation}\label{equality1}
\cig(G(q_0), -\nabla\Phi^n(\cdot,\lambda_\pm)) = \cig(G(q_0),-\nabla\Phi^{n_0}(\cdot,\lambda_\pm)),
\end{equation}
see \cite{IZYDOREK} for more details.

\nt
For simplicity of notation we set $H=G_{q_0}=\Gamma_{q_0}\times S^1$ and let $\ds \Psi^n_{\pm}=\Psi_{\pm|\bH^n}:\bH^n\rightarrow\bR,\, \Psi_{\pm}=\Phi(\cdot,\lambda_\pm)_{|\bH}:\bH\rightarrow\bR$ and $\ds \bH = T^\bot_{q_0}G(q_0), \,\bH^n=T_{q_0}^\bot \Gamma(q_0) \oplus \bigoplus_{k=1}^{n}\bH_k.$
Analogously like in Lemma 4.1 of \cite{PCHRYST}, using the same $H$-equivariant gradient homotopy, it follows that for any $n\geq n_0$ we have
\begin{equation} \label{equality2}
\cih(\{q_0\},-\nabla\Psi^n_{\pm})=
\cih(\{q_0\},-\nabla\Psi^{n_0}_{\pm})
\end{equation}
and 
$
\cih(\{q_0\},-\nabla\Psi^{n_0}_{\pm})=S^{\bH^+_{1,\pm}\oplus\bW^+}
$
where the spectral decomposition of $\bH^{n_0}$ for the isomorphism $-\nabla^2 \Psi^{n_0}_{\pm}(q_0)$ is the following $\ds \bH^{n_0}=\bH_1\oplus\left(T_{q_0}^\bot \Gamma(q_0) \oplus \bigoplus_{k=2}^{n_0} \bH_k \right) = (\bH^-_{1,\pm}\oplus\bH^+_{1,\pm})\oplus(\bW^-\oplus\bW^+)$
and $\dim \bH^+_{1,-} \neq \dim\bH^+_{1,+}.$

From now on we consider two cases: $\Gamma_{q_0}\approx T^{l_0}$ or $\Gamma$ is  abelian.

%%%%%%%%%%%%%%%%%%%%%%%%%%%%%%%%%%%%%%%%%%%%%%%%%%%%%%%%%%%%%%%%%%%%%%%%%%%%%%%%%%%%%%%%%%%%%%%%%%%%%%%%%%%%
%%%%%%%%%%%%%%%%%%%%%%%%%%%%%%%%%%%%%%%%%%%%%%%%%%%%%%%%%%%%%%%%%%%%%%%%%%%%%%%%%%%%%%%%%%%%%%%%%%%%%%%%%%%%

\subsubsection{CASE: $\Gamma_{q_0}\approx T^{l_0}$}
\numsubsubsec

For $l_0=0$ we have proved this theorem in \cite{PCHRYST}. Let $l_0 > 0.$ 
In this case $H \approx T^{l_0} \times S^1.$
Analogously like in Lemma 4.1 of \cite{PCHRYST} we have $\bH^+_{1,+}=\bH^+_{1,-}\oplus \bU$ where $\bU=\{a\cos t+b\sin t: a,b \in \bV_{\nabla^2U(q_0)}(\beta_{j_0}^2)\}$ and 
$\bV_{\nabla^2 U(q_0)}(\beta_{j_0}^2)$ is the eigenspace of $\nabla^2U(q_0)$ corresponding to the eigenvalue $\beta_{j_0}^2$.
Note that $\bV_{\nabla^2U(q_0)}(\beta_{j_0}^2)$ is an orthogonal representation of the torus  $T^{l_0}.$  
Since $\bU^H=\{0\},$ see  Lemma \ref{lem_fixedpoints},   $\bU$ is a nontrivial, fixed point free, orthogonal representation  
of  $H$.
Thus
\[\cih(\{q_0\},-\nabla\Psi^{n_0}_{-})=S^{\bH^+_{1,-}\oplus\bW^+}, \,
\cih(\{q_0\},-\nabla\Psi^{n_0}_{+})=S^{\bH^+_{1,-}\oplus \bU \oplus \bW^+}
\] and by Theorem \ref{theo_rozne} we obtain  $\chig(G^+\wedge_{H}S^{\bH^+_{1,-}\oplus\bW^+}) \neq \chig(G^+\wedge_{H}S^{\bH^+_{1,-} \oplus \bU \oplus \bW^+}).$ 
Hence by Theorem \ref{thmconley}.\eqref{thmconleyb}, the equalities \eqref{equality1} and \eqref{equality2} for any $n\geq n_0$ we obtain 
\begin{gather*}
\chig(\cig(G(q_0),-\nabla\Phi^{n}(\cdot,\lambda_-))) =\chig( G^+\wedge_{H}\cih(\{q_0\},-\nabla\Psi^{n}_{-}))\neq \\
\neq \chig( G^+\wedge_{H}\cih(\{q_0\},-\nabla\Psi^{n}_{+})) = \chig(\cig(G(q_0),-\nabla\Phi^{n}(\cdot,\lambda_+))) ,
\end{gather*}
and consequently 
$\cig(G(q_0),-\nabla\Phi^{n}(\cdot,\lambda_-)) \neq \cig(G(q_0),-\nabla\Phi^{n}(\cdot,\lambda_+))$ for any $n\geq n_0,$
which completes the proof of the first case.

%%%%%%%%%%%%%%%%%%%%%%%%%%%%%%%%%%%%%%%%%%%%%%%%%%%%%%%%%%%%%%%%%%%%%%%%%%%%%%%%%%%%%%%%%%%%%%%%%%%%%%%%%%%%
%%%%%%%%%%%%%%%%%%%%%%%%%%%%%%%%%%%%%%%%%%%%%%%%%%%%%%%%%%%%%%%%%%%%%%%%%%%%%%%%%%%%%%%%%%%%%%%%%%%%%%%%%%%%

\subsubsection{CASE: $\Gamma$ is abelian}
\numsubsubsec

By Example \ref{remarkamissible1}, the pair $(G,H)=(\Gamma \times \sone, \Gamma_{q_0} \times \sone)$ is admissible.
In view of  Theorem \ref{thmconley}.\eqref{thmconleyc} to prove this case it is enough to show that for any $n\geq n_0$
\[\chih(\cih(\{q_0\},-\nabla\Psi^n_{-})) \neq \chih(\cih(\{q_0\},-\nabla\Psi^n_{+})).\]
Let $i^\star:U(H)\rightarrow U(S^1)$  be the ring homomorphism induced by the inclusion $i:\{e\}\times S^1\rightarrow H,$ see Remark \ref{remarkringhom}.
Like in the proof of the previous case 
$\cih(\{q_0\},-\nabla\Psi^{n_0}_{-})=S^{\bH^+_{1,-}\oplus\bW^+}$ and 
$\cih(\{q_0\},-\nabla\Psi^{n_0}_{+})=S^{\bH^+_{1,-}\oplus \bU \oplus \bW^+}$
where $\bU$ is a nontrivial representation of $S^1\approx \{e\}\times S^1\subset H$ and the action of $S^1$ is given by shift in time.
Since $\bU$ is a nontrivial representation of $S^1$ and by the equality \eqref{equality2}, for any $n\geq n_0$ we obtain 
\begin{gather*}
i^\star(\chih(\cih(\{q_0\},-\nabla\Psi^{n}_{-})))=i^\star(\chih(S^{\bH^+_{1,-}\oplus\bW^+}))=\chi_{S^1}(S^{\bH^+_{1,-}\oplus\bW^+})\neq \\
\neq \chi_{S^1}(S^{\bH^+_{1,-}\oplus\bU\oplus\bW^+})= i^\star(\chih(S^{\bH^+_{1,-}\oplus\bU\oplus\bW^+})) = i^\star(\chih(\cih(\{q_0\},-\nabla\Psi^{n}_+))),
\end{gather*}
see Lemma 4.1 of \cite{PCHRYST}, and consequently 
$\chih(\cih(\{q_0\},-\nabla\Psi^{n}_{-}))\neq \chih(\cih(\{q_0\},-\nabla\Psi^{n}_+))$ for any $n\geq n_0,$  which completes the proof of the second case.

%%%%%%%%%%%%%%%%%%%%%%%%%%%%%%%%%%%%%%%%%%%%%%%%%%%%%%%%%%%%%%%%%%%%%%%%%%%%%%%%%%%%%%%%%%%%%%%%%%%%%%%%%%
%%%%%%%%%%%%%%%%%%%%%%%%%%%%%%%%%%%%%%%%%%%%%%%%%%%%%%%%%%%%%%%%%%%%%%%%%%%%%%%%%%%%%%%%%%%%%%%%%%%%%%%%%%%%
%%%%%%%%%%%%%%%%%%%%%%%%%%%%%%%%%%%%%%%%%%%%%%%%%%%%%%%%%%%%%%%%%%%%%%%%%%%%%%%%%%%%%%%%%%%%%%%%%%%%%%%%%%%%
%%%%%%%%%%%%%%%%%%%%%%%%%%%%%%%%%%%%%%%%%%%%%%%%%%%%%%%%%%%%%%%%%%%%%%%%%%%%%%%%%%%%%%%%%%%%%%%%%%%%%%%%%%%%
%%%%%%%%%%%%%%%%%%%%%%%%%%%%%%%%%%%%%%%%%%%%%%%%%%%%%%%%%%%%%%%%%%%%%%%%%%%%%%%%%%%%%%%%%%%%%%%%%%%%%%%%%%%%
%%%%%%%%%%%%%%%%%%%%%%%%%%%%%%%%%%%%%%%%%%%%%%%%%%%%%%%%%%%%%%%%%%%%%%%%%%%%%%%%%%%%%%%%%%%%%%%%%%%%%%%%%%%%
%%%%%%%%%%%%%%%%%%%%%%%%%%%%%%%%%%%%%%%%%%%%%%%%%%%%%%%%%%%%%%%%%%%%%%%%%%%%%%%%%%%%%%%%%%%%%%%%%%%%%%%%%%%%
%%%%%%%%%%%%%%%%%%%%%%%%%%%%%%%%%%%%%%%%%%%%%%%%%%%%%%%%%%%%%%%%%%%%%%%%%%%%%%%%%%%%%%%%%%%%%%%%%%%%%%%%%%%%
%%%%%%%%%%%%%%%%%%%%%%%%%%%%%%%%%%%%%%%%%%%%%%%%%%%%%%%%%%%%%%%%%%%%%%%%%%%%%%%%%%%%%%%%%%%%%%%%%%%%%%%%%%%%
%%%%%%%%%%%%%%%%%%%%%%%%%%%%%%%%%%%%%%%%%%%%%%%%%%%%%%%%%%%%%%%%%%%%%%%%%%%%%%%%%%%%%%%%%%%%%%%%%%%%%%%%%%%%

\subsection{Proof of Theorem \ref{main-theo2}}  \label{proofthm2}
\numsubsec
 
We follow the notation used in \cite{PCHRYST2}.
Fix $\beta_{j_0}$ satisfying the assumptions of Theorem \ref{main-theo2} and choose $\varepsilon>0$ such that $[\lambda_-,\lambda_+]\cap\Lambda=\{\frac{1}{\beta_{j_0}}\}$ where $\lpm=\frac{1\pm\varepsilon}{\beta_{j_0}}.$
To prove this theorem we have  to show  a change of the $G$-equivariant Conley index of the orbit $G(q_0)$ at level $\frac{1}{\beta_{j_0}},$ i.e.
$\sci_{G}\left(G(q_0),-\nabla\Phi(\cdot,\lambda_-)\right)\neq \sci_{G}(G(q_0),-\nabla\Phi(\cdot,\lambda_+)).$

In the proof we will apply the generalized equivariant Euler characteristic 
$\Upsilon_G : [GS(\xi)] \rightarrow U(G)$ defined for the category $[GS(\xi)]$ of $G$-homotopy types of $G$-spectra, see \cite{GORY1} for definition and properties. It is a  natural extension of the equivariant Euler characteristic defined for finite, pointed $G$-CW-complexes. 
To prove this theorem it is sufficient to show that $\Upsilon_G(\sci_G(G(q_0),-\nabla \Phi(\cdot,\lambda_-))) \neq  \Upsilon_G(\sci_G(G(q_0),-\nabla\Phi(\cdot,\lambda_+))).$

Let $H=G_{q_0}=\Gamma_{q_0} \times \sone$ and  $\h=T_{q_0}^{\perp} G(q_0) \subset \h^1_{2\pi}.$
Since
$\sci_{H}(\{q_0\},-\nabla\Psi_{\pm})=\sci_{H}(\{0\},-\nabla\tilde{\Psi}_{\pm})$
where
$\tilde{\Psi}_{\pm}(q)=\Psi_{\pm}(q+q_0)$ and $\Psi_\pm=\Phi(\cdot,\lambda_{\pm})_{\mid \bH} : \bH \to \bR,$
we study $\sci_H(\{0\},-\nabla\tilde{\Psi}_{\pm}).$
Set  $\bH= \ker \nabla^2 \tilde{\Psi}_\pm(0) \oplus \im \nabla^2\tilde{\Psi}_\pm(0) = Null \oplus Range$
where $Null$ is a finite-dimensional, orthogonal representation of $H$ such that $Null^{S^1}=Null,$ i.e.
$Null\subset\bH_0$ consists of constant functions.
Moreover, $Range$ is an infinite-dimensional, orthogonal representation of $H.$ Therefore if $\omega :Null\rightarrow Range$ is a $H$-equivariant map, then $\omega(Null)\subset Range^{S^1}=\bH_0\cap Range,$ i.e. the image
of $\omega$ also consists of constant functions. It is easy to prove Theorems 2.5.1 and 2.5.2 of \cite{PCHRYST2} for 
$H=G_{q_0}=\Gamma_{q_0} \times \sone$ and with the assumption (F.3) of \cite{PCHRYST2} replaced by $Null\subset\bH_0.$
To simplify notation, set $A_\pm=\nabla^2\tilde{\Psi}_\pm(0)_{|Range}.$ 
Applying this slight modification of Theorem 2.5.2 of \cite{PCHRYST2} we obtain $\varepsilon_0>0$ and  $H$-equivariant gradient homotopies
$\nabla \mathcal{H}_\pm: (B_{\varepsilon_0}(Null)\times B_{\varepsilon_0}(Range))\times[0,1]\rightarrow\bH $
satisfying the following conditions:
\begin{enumerate}
\item[(1)] $\nabla \mathcal{H}_\pm((u,v),t)=A_\pm v - \nabla\xi_\pm((u,v),t)$ for $t\in[0,1]$ where $\nabla \xi_\pm :\bH \times [0,1] \rightarrow \bH$ is compact and $H$-equivariant such that $\nabla\xi_\pm(0,t)=0, \nabla^2 \xi_\pm(0,t)=0$ for any $t\in[0,1],$ 
\item[(2)] $(\nabla \mathcal{H}_\pm)^{-1}(0)\cap (B_{\varepsilon_0}(Null)\times B_{\varepsilon_0}(Range))\times[0,1]=\{0\}\times[0,1],$ i.e. $0$ is an isolated critical point of $\nabla \mathcal{H}_\pm(\cdot,t)$ for any $t\in[0,1],$
\item[(3)] $\nabla \mathcal{H}_\pm((u,v),0)=\nabla\tilde{\Psi}_{\pm}(u,v),$
\item[(4)] there exist $H$-equivariant, gradient mappings $\nabla\varphi_\pm:B_{\varepsilon_0}(Null)\rightarrow Null$ such that 
$\nabla \mathcal{H}_\pm((u,v),1)=(\nabla\varphi_\pm(u), A_\pm v)$ for all $(u,v)\in B_{\varepsilon_0}(Null)\times B_{\varepsilon_0}(Range),$
\item[(5)] $\varphi_{\pm}(u)=\widetilde{\Psi}_{\pm}(u,w(u)), \nabla\varphi_{\lambda_{\pm}}(u)=P\nabla \widetilde{\Psi}_{\lambda_{\pm}} (u,w(u))$ is $H$-equivariant and $P:\bH\rightarrow Null$ is the $H$-equivariant, orthogonal projection 
and moreover, since $\omega(Null)\subset \bH_0,$
$\varphi_{\pm}(u) = -2\pi \lambda^2_{\pm} \tilde{U}(u,w(u))$ 
where $\tilde{U}: T^{\bot}_{q_0}\Gamma(q_0)\rightarrow\bR$ is defined by $\tilde{U}(q)=U(q+q_0).$
\end{enumerate}
Since the $H$-equivariant Conley index of $\{0\} \subset \bH$ is invariant
under the homotopy $\nabla \mathcal{H}_\pm(\cdot,t), 0 \leq t \leq 1$ and $\nabla \mathcal{H}_\pm(\cdot,1)$ is a product map, we obtain the following 
\begin{equation} \label{equal1}
\sci_{H}(\{0\},-\nabla \widetilde{\Psi}_{\pm})
=\sci_{H}(\{0\},-\nabla\varphi_{\pm})\wedge \sci_{H}(\{0\},-A_\pm) \in [HS(\xi)] .
\end{equation}

\nt
Since $Null\subset\bH_0$ is finite-dimensional, applying Remark 2.3.4 of \cite{PCHRYST2}, we have
\[\sci_{H}(\{0\},-\nabla\varphi_{\pm})=\cih(\{0\},-\nabla\varphi_{\pm}),\]
and so 
\begin{equation}\label{equal2}
\upsh(\sci_{H}(\{0\},-\nabla\varphi_{\pm}))=\chih(\cih(\{0\},-\nabla\varphi_{\pm})).
\end{equation}

\nt
From the property (5) it follows that $0 \in Null$ is an isolated critical point of $\varphi_\pm$ which is a local
maximum.
Therefore, by Lemma \ref{rabinowitzco}, we obtain that 
\begin{equation} \label{equal3}
\chih(\cih(\{0\},-\nabla\varphi_{\pm}))=\chih(S^{Null})\in U(H).
\end{equation}
Applying the generalized equivariant Euler characteristic $\upsh$ to the equality \eqref{equal1} and combining the equalities \eqref{equal2} and \eqref{equal3}  we get
\begin{multline} \label{inequality1}
\upsh\left(\sci_{H}(\{q_0\},-\nabla \Psi_\pm)\right) = \upsh \left(\sci_{H}(\{0\},-\nabla \widetilde{\Psi}_\pm)\right)=\\
=\upsh( \sci_{H}(\{0\},-\nabla \varphi_\pm))\ast\upsh\left(\sci_{H}(\{0\},-A_\pm)\right) 
=\chih(S^{Null})\ast \upsh\left(\sci_{H}(\{q_0\},-\nabla\Pi_{\pm})\right)
\end{multline}
where
$\Pi_{\pm}:Range \to\bR$ is given by $\Pi_{\pm}(q) = \frac{1}{2}\langle \nabla^2 \Psi_\pm(q_0)_{\mid Range}(q-q_0), q-q_0 \rangle_{\h^1_{2\pi}} .$ 
Let 
$\Pi_{\pm}^n = \Pi_{\pm|\bH^n} : \bH^n\rightarrow\bR$ and $\ds \bH^n=(T_{q_0}^\bot\Gamma(q_0) \oplus \bigoplus_{k=1}^{n}\bH_k)\cap Range
=\left(T_{q_0}^{\perp} \Gamma(q_0)\ominus Null \right) \oplus \bigoplus_{k=1}^{n} \bH_k\subset Range \subset \bH.$
Analogously like in Lemma 3.3.2 of \cite{PCHRYST2} it follows that there exists $n_0\in\bN$ such that for any $n\geq n_0$
\begin{equation}\label{equality4}
\cih(\{q_0\},-\nabla\Pi^n_{\pm})=\cih(\{q_0\},-\nabla\Pi^{n_0}_{\pm})
\end{equation}
and
\beq \label{inequality2}
\cih(\{q_0\},-\nabla\Pi^{n_0}_{\pm})= S^{\bH^+_{1,\pm}\oplus\bW^+}
\eeq
where the spectral decomposition of $\bH^{n_0}$ given by the isomorphism $-\nabla\Pi^{n_0}_{\lambda_\pm}$ is the following \[\bH^{n_0}=\bH_1\oplus\left((T_{q_0}^\bot\Gamma(q_0)\ominus Null)\oplus\bigoplus_{k=2}^{n_0}\bH_k\right)=(\bH^-_{1,\pm}\oplus\bH^+_{1,\pm})\oplus(\bW^-\oplus\bW^+)\]
and $\dim \bH^+_{1,-} \neq \dim\bH^+_{1,+}.$
Therefore the $H$-equivariant Conley indexes $\sci_H(\{q_0\},-\nabla\Pi_{-})$ and $\sci_H(\{q_0\},-\nabla\Pi_{\lambda_+})$ are the $H$-homotopy types of $H$-spectra of the same type $\xi=(\bV_n)_{n=0}^{\infty}$. 
We put $\bV^n=\bV_0 \oplus \bV_1 \oplus \cdots \oplus \bV_n.$ 
Now basing on the definition of the generalized equivariant Euler characteristic 
$\upsh : [HS(\xi)] \rightarrow U(H)$ we obtain that 
\beq \label{inequality11}
\upsh(\sci_H(\{q_0\},-\nabla\Pi_{\pm}))=\chih\left(S^{\bV^{n_0-1}}\right)^{-1}\star \chih(\cih(\{q_0\},-\nabla\Pi^{n_0}_{\pm})).
\eeq
Recall that $\chih(S^{Null})$ is an invertible element of $U(H),$ see Theorem 3.5 of \cite{GERY}. 
Since $Null$ does not depend on the levels $\lambda_-$ and $\lambda_+,$ combining the equalities \eqref{inequality1} and \eqref{inequality11}, we have that 
the condition 
$\upsh\left(\sci_{H}(\{q_0\},-\nabla \Psi_-)\right)\neq \upsh\left(\sci_{H}(\{q_0\},-\nabla\Psi_+)\right)$ is equivalent to the condition
$\cih(\{q_0\},-\nabla\Pi_{-}^n)\neq \cih(\{q_0\},-\nabla \Pi_{+}^n)$ for any $n\geq n_0.$ 

From now on we consider two cases: $\Gamma_{q_0}\approx T^{l_0}$ or $\Gamma$ is abelian.

%%%%%%%%%%%%%%%%%%%%%%%%%%%%%%%%%%%%%%%%%%%%%%%%%%%%%%%%%%%%%%%%%%%%%%%%%%%%%%%%%%%%%%%%%%%%%%%%%%%%%%%%%%
%%%%%%%%%%%%%%%%%%%%%%%%%%%%%%%%%%%%%%%%%%%%%%%%%%%%%%%%%%%%%%%%%%%%%%%%%%%%%%%%%%%%%%%%%%%%%%%%%%%%%%%%%%%%

\subsubsection{CASE: $\Gamma_{q_0}\approx T^{l_0}$.}
For $l_0=0$ we have proved this theorem in \cite{PCHRYST2}. Let $l_0 > 0.$ 
In this case $H \approx T^{l_0} \times S^1.$
Analogously like in Lemma 3.3.2 of \cite{PCHRYST2} we have $\bH^+_{1,+}=\bH^+_{1,-}\oplus \bU$ where $\bU=\{a\cos t+b\sin t:  a,b \in \bV_{\nabla^2U(q_0)}(\beta_{j_0}^2)\}$ and  
$\bV_{\nabla^2 U(q_0)}(\beta_{j_0}^2)$ is the eigenspace of $\nabla^2U(q_0)$ corresponding to the eigenvalue $\beta_{j_0}^2$.
From  Lemma \ref{lem_fixedpoints} we get $\bU^H=\{0\},$ and consequently that $\bU$ is a nontrivial representation of $H$. Thus
\begin{equation} 
\cih(\{q_0\},-\nabla\Pi^{n_0}_{-})=S^{\bH^+_{1,-}\oplus\bW^+}, \,
\cih(\{q_0\},-\nabla\Pi^{n_0}_{+})=S^{\bH^+_{1,-}\oplus\bU \oplus \bW^+}.
\end{equation} 
By the equalities \eqref{equal1} and \eqref{equality4} we obtain for any $n\geq n_0$ 
\[\cih(\{q_0\},-\nabla\Psi^{n}_\pm)=S^{Null}\wedge \cih(\{q_0\},-\nabla\Pi_\pm^{n}).\]
According to Theorem 2.4.2 of \cite{PCHRYST2} for isolated critical orbits 
we have 
\[\cig(G(q_0),-\nabla\Phi^{n}(\cdot,\lambda_\pm)) = G^+\wedge_{H}\cih(\{q_0\},-\nabla\Psi^{n}_{\pm})\] 
for any $n\geq n_0.$
Therefore since $\bU$ is a nontrivial representation of $H$, applying  Theorem \ref{theo_rozne},  we get for any $n\geq n_0$ 
\begin{gather*}\chig(\cig(G(q_0),-\nabla\Phi^{n}(\cdot,\lambda_-))) = \chig(G^+ \wedge_{H} S^{Null\oplus\bH^+_{1,-} \oplus \bW^+}) \neq \\
\chig(G^+\wedge_{H} S^{Null\oplus\bH^+_{1,-}\oplus\bU\oplus\bW^+}) = \chig(\cig(G(q_0),-\nabla\Phi^{n}(\cdot,\lambda_+))).
\end{gather*}
Consequently $\cig(G(q_0),-\nabla\Phi^{n}(\cdot,\lambda_-)) \neq \cig(G(q_0),-\nabla\Phi^{n}(\cdot,\lambda_+))$ for any $n\geq n_0$
which implies that  
$
\sci_{G}\left(G(q_0),-\nabla\Phi(\cdot,\lambda_-)\right)\neq \sci_{G}(G(q_0),-\nabla\Phi(\cdot,\lambda_+)).
$ This completes the proof of the first case.

%%%%%%%%%%%%%%%%%%%%%%%%%%%%%%%%%%%%%%%%%%%%%%%%%%%%%%%%%%%%%%%%%%%%%%%%%%%%%%%%%%%%%%%%%%%%%%%%%%%%%%%%%%%%
%%%%%%%%%%%%%%%%%%%%%%%%%%%%%%%%%%%%%%%%%%%%%%%%%%%%%%%%%%%%%%%%%%%%%%%%%%%%%%%%%%%%%%%%%%%%%%%%%%%%%%%%%%%%

\subsubsection{CASE: $\Gamma$ is abelian}
\numsubsubsec
By Example \ref{remarkamissible1}, the pair $(G,H)=(\Gamma \times \sone, \Gamma_{q_0} \times \sone)$ is admissible. To prove this case it is sufficient to show that
\[\Upsilon_{H}(\sci_{H}(\{q_0\},-\nabla\Psi_{-}))\neq\Upsilon_{H}(\sci_{H}(\{q_0\},-\nabla\Psi_{\lambda_+})),\]
see Theorem 2.4.3 of \cite{PCHRYST2}. 
We have proved that this condition is equivalent  to
\[\cih(\{q_0\},-\nabla\Pi_{-}^n)\neq \cih(\{q_0\},-\nabla\Pi_{+}^n)\] for any $n\geq n_0.$
Let $i^\star:U(H) \rightarrow U(\sone)$ be the ring homomorphism induced by the inclusion  $i:\{e\}\times \sone \rightarrow H,$ see  Remark \ref{remarkringhom}.
Like in the proof of the previous case we have
$\cih(\{q_0\},-\nabla\Pi_{-}^{n_0}) =  $ $S^{\bH^+_{1,-}\oplus\bW^+}$ and $\cih(\{q_0\},-\nabla\Pi_{-}^{n_0})=S^{\bH^+_{1,-}\oplus\bU\oplus\bW^+}$ where $\bU$ is a nontrivial representation of $S^1.$
Since $U$ is a nontrivial representation of $S^1$ and by the equality \eqref{equality4}, for any $n\geq n_0$ we obtain the following 
\begin{gather*}
i^\star(\chih(\cih(\{q_0\},-\nabla\Pi^{n}_{-})))=i^\star(\chih(S^{\bH^+_{1,-} \oplus \bW^+})) = \chi_{S^1}(S^{\bH^+_{1,-}\oplus\bW^+})\neq \\
\neq \chi_{S^1}(S^{\bH^+_{1,-}\oplus\bU\oplus\bW^+})= i^\star(\chih(S^{\bH^+_{1,-}\oplus\bU\oplus\bW^+}))=i^\star(\chih(\cih(\{q_0\},-\nabla\Pi^{n}_{\lambda_+}))),
\end{gather*}
see Lemma 3.3.2 of \cite{PCHRYST2}.  Consequently 
$\cih(\{q_0\},-\nabla\Pi^{n}_{-})\neq \cih(\{q_0\},-\nabla\Pi^{n}_{+})$ for any $n \geq n_0,$
which completes the proof of the second case.

%%%%%%%%%%%%%%%%%%%%%%%%%%%%%%%%%%%%%%%%%%%%%%%%%%%%%%%%%%%%%%%%%%%%%%%%%%%%%%%%%%%%%%%%%%%%%%%%%%%%%%%%%%
%%%%%%%%%%%%%%%%%%%%%%%%%%%%%%%%%%%%%%%%%%%%%%%%%%%%%%%%%%%%%%%%%%%%%%%%%%%%%%%%%%%%%%%%%%%%%%%%%%%%%%%%%%%%
%%%%%%%%%%%%%%%%%%%%%%%%%%%%%%%%%%%%%%%%%%%%%%%%%%%%%%%%%%%%%%%%%%%%%%%%%%%%%%%%%%%%%%%%%%%%%%%%%%%%%%%%%%%%
%%%%%%%%%%%%%%%%%%%%%%%%%%%%%%%%%%%%%%%%%%%%%%%%%%%%%%%%%%%%%%%%%%%%%%%%%%%%%%%%%%%%%%%%%%%%%%%%%%%%%%%%%%%%
%%%%%%%%%%%%%%%%%%%%%%%%%%%%%%%%%%%%%%%%%%%%%%%%%%%%%%%%%%%%%%%%%%%%%%%%%%%%%%%%%%%%%%%%%%%%%%%%%%%%%%%%%%%%
%%%%%%%%%%%%%%%%%%%%%%%%%%%%%%%%%%%%%%%%%%%%%%%%%%%%%%%%%%%%%%%%%%%%%%%%%%%%%%%%%%%%%%%%%%%%%%%%%%%%%%%%%%%%
%%%%%%%%%%%%%%%%%%%%%%%%%%%%%%%%%%%%%%%%%%%%%%%%%%%%%%%%%%%%%%%%%%%%%%%%%%%%%%%%%%%%%%%%%%%%%%%%%%%%%%%%%%%%
%%%%%%%%%%%%%%%%%%%%%%%%%%%%%%%%%%%%%%%%%%%%%%%%%%%%%%%%%%%%%%%%%%%%%%%%%%%%%%%%%%%%%%%%%%%%%%%%%%%%%%%%%%%%
%%%%%%%%%%%%%%%%%%%%%%%%%%%%%%%%%%%%%%%%%%%%%%%%%%%%%%%%%%%%%%%%%%%%%%%%%%%%%%%%%%%%%%%%%%%%%%%%%%%%%%%%%%%%
%%%%%%%%%%%%%%%%%%%%%%%%%%%%%%%%%%%%%%%%%%%%%%%%%%%%%%%%%%%%%%%%%%%%%%%%%%%%%%%%%%%%%%%%%%%%%%%%%%%%%%%%%%%%


\begin{thebibliography}{99}

\bibitem{BERGER1} M. S. Berger. Bifurcation theory and the type numbers of Marston Morse. Proc. Nat. Acad. Sci. U.S.A.,
69:1737-1738, 1972.

\bibitem{BERGER2}  M. S. Berger. Nonlinearity and functional analysis: Lectures on nonlinear problems in mathematical analysis.
Pure and Applied Mathematics. Academic Press, New York-London, 1977.

\bibitem{UGORY} J. Fura, A. Ratajczak, and S. Rybicki. Existence and continuation of periodic solutions of autonomous
Newtonian systems. J. Differential Equations, 218(1):216-252, 2005.

\bibitem{GEBA} K. G\k{e}ba. Degree for gradient equivariant maps and equivariant Conley index. In M. Matzeu and A. Vignoli,
editors, Topological Nonlinear Analysis II: Degree, Singularity and Variations, volume 27 of Progr. Nonlinear
Differential Equations Appl., pages 247-272. Birkh{\"a}user Boston, 1997.

\bibitem{GERY} K. G\k{e}ba and S. Rybicki. Some remarks on the Euler ring U(G). J. Fixed Point Theory Appl., 3(1):143-158,
2008.

\bibitem{GORY1} A. Go\l{}\k{e}biewska and S. Rybicki. Equivariant Conley index versus degree for equivariant gradient maps. Discrete
Contin. Dyn. Syst. Ser. S, 6(4):985-997, 2013.

\bibitem{HENRARD}  J. Henrard. Lyapunov's center theorem for resonant equilibrium. J. Differential Equations, 14:431-441, 1973.

\bibitem{IZYDOREK} M. Izydorek. Equivariant Conley index in Hilbert spaces and applications to strongly indefinite problems.
Nonlinear Anal., 51(1):33-66, 2002.

\bibitem{KBO} K. Kawakubo. The theory of transformation groups. The Clarendon Press, Oxford University Press, New
York, 1991.

\bibitem{GLG} G. L{\'o}pez Garza and S. Rybicki. Equivariant bifurcation index. Nonlinear Anal., 73(9):2779-2791, 2010.

\bibitem{MAWI} J. Mawhin and M. Willem. Critical point theory and Hamiltonian systems, volume 74 of Applied Mathematical
Sciences. Springer-Verlag, New York, 1989.

\bibitem{PCHRYST} E. P{\'e}rez-Chavela, S. Rybicki, and D. Strzelecki. Symmetric Liapunov center theorem. Calc. Var. Partial
Differential Equations, 56(2):26, 2017.

\bibitem{PCHRYST2} E. P{\'e}rez-Chavela, S. Rybicki, and D. Strzelecki. Symmetric Liapunov center theorem for minimal orbit. J.
Differential Equations, 265(3):752-778, 2018.

\bibitem{RAB4} P. H. Rabinowitz. A note on topological degree for potential operators. J. Math. Anal. Appl., 51(2):483-492,
1975.

\bibitem{DIECK0} T. tom Dieck. Transformation groups and representation theory, volume 766 of Lecture Notes in Mathematics.
Springer-Verlag, Berlin-New York, 1979.

\bibitem{DIECK} T. tom Dieck. Transformation groups, volume 8 of De Gruyter Studies in Mathematics. Walter de Gruyter
and Co., Berlin, 1987.

\end{thebibliography}
\end{document}